\PassOptionsToPackage{usenames,dvipsnames}{xcolor}
\documentclass{infor}
\usepackage{graphicx}
\usepackage{booktabs}
\usepackage{url}
\usepackage{amsmath}
\usepackage{amssymb}
\usepackage{placeins}
\usepackage{algorithm}
\usepackage[noend]{algpseudocode}
\usepackage{subcaption}
\usepackage{mathtools}

\DeclareMathOperator{\argmax}{\mathrm{argmax}}

\volume{0}
\issue{0}
\pubyear{20XX}
\articletype{research-article}
\doi{0000}

\theoremstyle{remark}

\begin{document}
\begin{frontmatter}
\pretitle{Bi-Level Optimization to Enhance Intensity Modulated Radiation Therapy Planning}

\title{Bi-Level Optimization to Enhance Intensity Modulated Radiation Therapy Planning}

\author[a]{\inits{J.J}\fnms{Juan José} \snm{Moreno}\thanksref{c1}\ead[label=e1]{juanjomoreno@ual.es}\bio{bio1}}
\author[a,b]{\inits{S.}\fnms{Sav\'ins} \snm{Puertas-Mart\'in}\ead[label=e2]{savinspm@ual.es}\bio{bio2}}
\author[a]{\inits{J.L.}\fnms{Juana L.} \snm{Redondo}\ead[label=e3]{jlredondo@ual.es}\bio{bio3}}
\author[a]{\inits{P.M.}\fnms{Pilar M.} \snm{Ortigosa}\ead[label=e4]{ortigosa@ual.es}\bio{bio4}}
\author[c]{\inits{A.}\fnms{Anna} \snm{Zawadzka}\ead[label=e5]{anna.zawadzka.fizyk@nio.gov.pl}\bio{bio5}}
\author[c]{\inits{P.}\fnms{Pawel} \snm{Kukolowicz}\ead[label=e6]{Pawel.Kukolowicz@nio.gov.pl}\bio{bio6}}
\author[d]{\inits{R.}\fnms{Robert} \snm{Szmurło}\ead[label=e7]{robert.szmurlo@pw.edu.pl}\bio{bio7}}
\author[e]{\inits{I.}\fnms{Ignacy} \snm{Kaliszewski}\ead[label=e8]{ignacy.kaliszewski@ibspan.waw.pl}\bio{bio8}}
\author[e]{\inits{J.}\fnms{Janusz} \snm{Miroforidis}\ead[label=e9] {janusz.miroforidis@ibspan.waw.pl}\bio{bio9}}
\author[a]{\inits{E.M.}\fnms{Ester M.} \snm{Garzón}\ead[label=e10]{gmartin@ual.es}\bio{bio10}}
\thankstext[type=corresp,id=c1]{Corresponding author.}
\address[a]{Dept. of Informatics, \institution{University of Almería}, \cny{Spain}}
\address[b]{Information School, \institution{University of Sheffield}, \cny{United Kingdom}}
\address[c]{Department of Medicine Physics, \institution{Maria Skłodowska Curie Memorial Cancer Centre and Institute of Oncology}, \cny{Poland}}
\address[d]{Faculty of Electrical Engineering, \institution{Warsaw University of Technology}, \cny{Poland}}
\address[e]{Systems Research Institute, \institution{Polish Academy of Sciences}, \cny{Poland}}


\begin{abstract}
Intensity Modulated Radiation Therapy is an effective cancer treatment. Models based on the Generalized Equivalent Uniform Dose (gEUD) provide radiation plans with excellent planning target volume coverage and low radiation for organs at risk. However, manual adjustment of the parameters involved in gEUD is required to ensure that the plans meet patient-specific physical restrictions. This paper proposes a radiotherapy planning methodology based on bi-level optimization. We evaluated the proposed scheme in a real patient and compared the resulting irradiation plans with those prepared by clinical planners in hospital devices. The results in terms of efficiency and effectiveness are promising.
\end{abstract}


\end{frontmatter}

\section{Introduction}

Intensity Modulated Radiation Therapy (IMRT) is an effective radiation therapy technique in cancer treatment. It requires personalized irradiation plans (RT plans) that specify 3D dose distributions which effectively destroy cancer cells while minimizing side effects on  healthy tissues. Subsequently, linear accelerators equipped with multileaf collimators deliver radiation beams to patients from several fixed angles. Each beam is decomposed into a regular grid of (thousands of) beamlets with varying intensities that can be individually controlled. Therefore, every RT plan is defined by the specific intensities of all the beamlets over all beams. When preparing an RT plan, the goal is twofold.
On the one hand, plans are sought that ensure deposing the prescribed doses in planning target volumes (PTVs), and, on the other hand, such plans should protect organs at risk (OARs). These two goals are contradictory and, as such, must be traded-off, which means that high quality RT plans, i.e., effective and at the same time causing as little harmful side effects as possible, require significant efforts from the RT planners \citep{Breedveld_et_all_2019}. This contradiction leads the way to multi-objective optimization in IMRT.

To obtain effective RT plans, various multi-objective optimization models with several radiation effect measures (or criteria, at least one measure for each PTV and OAR), physical or biological, have been proposed \citep{Breedveld_et_all_2019,Cho18,EHR10}. Physical measures capture characteristics of dose distributions. Usually, in PTVs, physical measures take the form of sums of absolute or squared values of the difference between the delivered and prescribed doses. On the other hand, biological measures refer to (estimated) biological effects of radiation in organs \citep{Olaf05}. Doses delivered to OARs are controlled by imposing physical constraints on acceptable average or maximal radiation doses per voxel (a voxel being a 3D box constituting the irradiation planning mesh).
In PTVs and OARs, the appropriateness of dose distributions is also assessed by the shapes of dose-volume histograms (DVHs). Favorable DVH shapes are enforced by appropriate dose-volume constraints, which as a rule rend the optimization problems nonconvex \citep{BORT99}.
Given the multiplicity of beams, beamlets, and voxels, the resulting models are computationally complex. The most handleable are linear or linearized models \citep{Rom06,Olaf06}. To handle nonconvex constraints, in \cite{Breedveld_Heimen_2017} the interior point method has been adapted. Another approach to cope with nonconvexity is to resort to heuristics, or a combination of heuristic and exact methods \citep{Lan12}. Stochastic approaches are an alternative \citep{AHMAD10,Moreno21}. Finally, convex formulations of dose-volume histograms constraints using general formulations of dose-volume constraints is another option \citep{FU2019,Romeijn2004AUF}.

Nowadays, capitalizing on the progress in oncology research, biologically-motivated measures of radiation effects play an important role in clinical practice. Despite the difficulties in their implementation for clinical use, they are recommended by the ICRU\footnote{\mbox{ }International Committee for Radiological Units.} Reports as an advanced level of treatment reporting, taking into account developing concepts. Measures of dose distribution effectiveness, such as the tumor control probability (TCP) and the normal tissue complication probability (NTCP) \citep{Alb06}, the biologically effective dose (BED) \citep{Sab15}, and the generalized equivalent uniform dose (gEUD) \citep{Breedveld_et_all_2019,NIEM97,Wu_02}, are gaining popularity. The gEUD metric, based on the linear-quadratic cell survival model, has been reported to be the most relevant for radiotherapy \citep{NIEM97,Wu_02,Olaf05}. Therefore, it is the focus of this work.

For gEUD-based RT planning optimization, optimal solutions can be efficiently computed by gradient methods, as proposed in \cite{Choi_2002}. However, for every PTV and OAR several parameters need to be set by the RT planner. Usually, the planner starts from a set of parameters recommended in the literature and later adjusts them to the particular patient case by trial and error, to satisfy the imposed physical constraints. The quality of the final RT plan depends heavily on the knowledge and skills of the planner.

Despite differences between physical and biological measures, both are concurrently used in clinical settings. We propose an RT planning methodology based on a bi-level optimization scheme where: on the lower level, a gEUD-based objective function with all its parameters fixed (including the gEUD parameters) is optimized by an (exact) gradient algorithm and, on the upper level, those parameters are optimized by an evolutionary algorithm.

The output of the proposed scheme is a collection of non-dominated RT plans. With varying priorities, several RT plans (solutions to the multi-objective optimization problem) are derived. Furthermore, to facilitate the analysis of these RT plans, we have also developed a decision tool that assists the planner in selecting the final plan.

\section{Computational methods and theory: IMRT planning based on the gEUD and physical constraints}

This section describes the proposed methodology for obtaining IMRT plans by combining gEUD and physical constraints based on oncology physicians' recommendations and individual patient anatomy.

To facilitate the reading, the notation used throughout the paper is summarized in the Table \ref{t1:notation}.

\begin{table}[ht]
	\centering
	\caption{Notation and naming conventions}
	\label{t1:notation}
	\setlength{\tabcolsep}{6pt}
	\begin{tabular}{ll}
	\toprule
	\textbf{Notation} & \textbf{Meaning} \\
		\midrule
		$x$ & fluence map\\
	    $D$ & deposition matrix (translates $x$ to doses in voxels) \\
		$D_j$ & $j$-th row of $D$\\
		$d(x)=Dx$   & vector of doses deposited in voxels \\
		$d_j(x)=D_j  x$   & dose deposited in voxel $j$ \\
		$T= \{ t \}$ & set of indices of PTVs \\
		$M_t$ & set of voxels in $t$-th PTV \\
		$EUD^0_t$ & prescribed uniform dose for $t$-th PTV \\
		
		$R= \{ r \}$ & set of indices of OARs\\
		$M_r$ & set of voxels in $r$-th OAR \\
		$EUD^0_r$ & maximal uniform dose for $r$-th OAR\\
		$P= \{ p \} \subseteq R$ & subset of high-priority protected OARs\\
		\bottomrule
	\end{tabular}
\end{table}

\subsection{gEUD-based IMRT planning} \label{sec:EUD_model}

RT planning can be viewed as an optimization process, in which the fluence maps, represented by vectors of non-negative numbers $x$, define the radiation intensities of individual beamlets. Deposition matrices $D$ translate fluencies $x$ to doses deposited in voxels, so that the doses in voxels are computed as product $Dx$. The RT planning goal is to compute fluencies $x$ that deposit prescribed and homogeneous doses to PTVs and acceptable doses to OARs.

gEUD is a biology-motivated measure to evaluate radiation effects, based on the concept of the uniform radiation dose delivered to a patient organ, that causes the same effect as a nonuniform dose \citep{NIEM97,Wu_02}. 
In the case of gEUD, radiation effects in a PTV or an OAR, both referred to as structure $s$, are evaluated by the following function that aggregates these effects over all voxels belonging to the structure $s$: 
	\begin{equation}
	\label{eq_EUDparam}
		gEUD_s(x, a_s) = \left( \frac{1}{|M_{s}|}\sum_{j \in M_s}d_j(x)^{a_s} \right)^{\frac{1}{a_s}}
	\end{equation}
where $|M_{s}|$ is the number of voxels of the structure $s$; $d_j(x)=D_j x$ is the radiation dose deposited in voxel $j$ of structure $s$ by fluence map $x$ and the parameter $a_s$ represents the radiation effect on the structure. Its value can be taken from the literature or it can be adjusted individually by trial and error for each patient case.

According to \cite{Wu_02}, clinically meaningful RT plans can be obtained by computing the maximum of the following function $F(x, \phi)$, built over the gEUD:
\begin{equation}
\label{EUDmodel}
\begin{array}{ll}
F(x, \phi)=&\prod_{t \in T}  \frac{1}{1 +\left( \frac{EUD_t^0}{gEUD_t(x,a_t)} \right)^{n_t}} 
\cdot \prod_{r \in R} \frac{1}{1 + \left( \frac{gEUD_r(x,a_r)}{EUD_r^0}\right)^{n_r}} 
\end{array}
\end{equation}
where $EUD_t^0$ is the prescribed dose for $t$-th PTV, ${EUD_r^{0}}$ is the maximum uniform dose at $r$-th OAR; $n_r$ and $n_t$ express the importance of the prescriptions for the corresponding structure; $\phi$ represents the set of parameters involved in the $F$ definition, i.e., $\phi$ is an instance of parameters $n_t, n_r, a_t, a_r$ and $EUD_r^0$ with $ t \in T, r \in R$.

As is, the objective function $F(x, \phi)$ only controls underdosage inside PTVs. However, overdosage control in those volumes is also important. For this purpose it is common to define complementary structures, introduced as "virtual PTVs" on \cite{Wu_02}, that are treated as OARs. In this work, for each PTV we have defined a virtual PTV. To lighten optimization costs, we have interrelated the respective parameters, namely, for each PTV $t$, the corresponding virtual PTV $r$ has $EUD^0_r = EUD^0_t + 1$, $a_r = -a_t$ and $n_r = n_t$.

In \cite{Choi_2002} convexity of function (\ref{eq_EUDparam}) was studied. It was shown that for a range of parameters, function (\ref{eq_EUDparam}) is a convex function of fluence maps ($x$). In \cite{Romeijn2004AUF} the convexity of (\ref{EUDmodel}) was analyzed and the conclusion was that it is convex for the range of values of the parameters $a_s$ considered in practice. This means that IMRT RT plans can be efficiently obtained by gradient methods, as suggested in \cite{Wu_02}.

RT planning based on maximization of function $F(x,\phi)$ goes as follows. At the start, the values of $n_t, n_r, a_t, a_r, \ t \in T, \ r \in R$, are selected according to suggestions from the literature. 
Next, in successive planning cycles, the parameters $a_p$ and $EUD_p^0, \ p \in P$, are manually tuned by the RT planner. This is a tiresome and time-consuming process, requiring high-expertise planners. In this work, we propose to select the gEUD parameters automatically by multi-objective optimization. In the next section, we present this idea in more detail.

\subsection{Multi-objective optimization for gEUD-based IMRT planning} \label{sec:Multiobjetive}

Automated parameter tuning can result in several, even many, RT plans generated in one planning session. We present a multi-objective optimization model which serves that purpose. We also outline a novel approach to solve this model by a hybrid optimization scheme, coupling exact and heuristic (evolutionary) optimization methods.

\subsubsection{The multi-objective optimization model}
\label{ssec:multiobjectivemodel}

Here we present a multi-objective optimization model capable of handling physical constraints imposed on RT plans. Such constraints are individually selected in the optimization process based on the individual patient's anatomy.

For structure $s$ we consider the following dose distribution statistical measures:
\begin{equation}
\begin{array}{l}
D_s^{min}(x)  =  \min_{j \in M_{s}} d_j(x)\\ \\
\overline {D_s(x)}  =  \frac{1}{|M_{s}|}\sum_{j \in M_s}d_j(x)  \\ \\
D_s^{max}(x)  = \max_{j \in M_{s}} d_j(x) \\ \\
\end{array}
\end{equation}

With these measures we define constraints in structures. Four constraints are defined for each PTV: 

\begin{equation}
\begin{array}{l}
D_t^{min}(x) \geq LB_t \\ \\
\overline {D_t(x)} \geq \overline {LB_t}\\ \\
\overline {D_t(x)} \leq \overline {UB_t}\\ \\
D_t^{max}(x) \leq UB_t
\end{array}
\end{equation}
where, for given $t$, $LB_t$ and $UB_t$ are the lower and upper bound for the dose in any voxel of the structure, $\overline {LB_t}, \overline {UB_t}$ are lower and upper bound for the average dose in the structure. 

Doses in parallel \footnote{\mbox{ }An organ is called {\it parallel} if its functionality is preserved despite partial radiation damage, e.g., the salivary gland.} OARs are constrained by upper bounds on the average dose in the structure: 
\begin{equation}
\begin{array}{l}
\overline {D_r(x)} \leq \overline {UB_r}\\ \\
\end{array}
\end{equation}
whereas doses in serial\footnote{\mbox{ }An organ is called {\it serial} if it loses its functionality completely if any part of it is damaged, e.g., the spinal cord.} OARs are constrained by upper bounds on the maximal dose in individual voxels of the structure:
\begin{equation}
\begin{array}{l}
D_r^{max}(x) \leq UB_r\\ \\
\end{array}
\end{equation}

Constraint violations in structures, PTVs and OARs, are captured by the following constraint violation functions:

\begin{equation}
\begin{array}{ll}
\vspace{1em}
C_s^{min}(x) =&
\begin{dcases}
    LB_s - D_s^{min}(x) & \text{if } LB_s \text{ is defined in } s \text{ and } LB_s > D_s^{min}(x)\\
    0              & \text{otherwise}
\end{dcases}\\
\vspace{1em}
\overline{C_s^{min}}(x) =&
\begin{dcases}
    \overline {LB_s} - \overline {D_s(x)} & \phantom{123}  \text{if } \overline {LB_s} \text{ is defined in } s \text{ and } \overline{LB_s} > \overline {D_s(x)}\\
    0              & \phantom{123}  \text{otherwise}
\end{dcases}\\
\vspace{1em}
\overline{C_s^{max}}(x) =&
\begin{dcases}
    \overline {D_s(x)} - \overline {UB_s} & \phantom{123} \text{if } \overline {UB_s} \text{ is defined in } s \text{ and } \overline{UB_s} < \overline {D_s(x)}\\
    0              & \phantom{123}  \text{otherwise}
\end{dcases}\\
\vspace{1em}
C_s^{max}(x) =&
\begin{dcases}
    D_s^{max}(x) - UB_s  & \text{if } UB_s \text{ is defined in } s \text{ and } UB_s < D_s^{max}(x)\\
    0              & \text{otherwise}
\end{dcases}\\
\end{array}
\end{equation}

Since obviously constraint violations should be minimized, we formulate the following objective functions. 
The first function, denoted $f_0$, aggregates constraint violations over all structures:
\begin{equation}
	\label{f0}
f_0(x) =  \sum_{s \in S} C_s^{min}(x) + \overline{C_s^{min}}(x) + \overline{C_s^{max}}(x) + C_s^{max}(x)
\end{equation}

The priorities to protect the distinguished subsets of OARs, indexed by $p, \ p \in P \subseteq R$, are represented in the model by separate objective functions being the averages of the deposited doses, $\overline{D_p(x)}$: 
\begin{equation}
	\label{fp1}
f_p(x) = \overline {D_p(x)}, \ \ p \in P,
\end{equation}
or
\begin{equation}
	\label{fp2}
	f_p(x) = D_p^{max}(x), \ \ p \in P.
\end{equation}

Clearly, low values of individual functions $f_0, f_1, \dots, f_{|P|}$ signal acceptable plans. 

The multi-objective optimization problem consists of $|P| + 1$ objective functions $f_0, \dots, f_{|P|}$, all of them to be minimized:
\begin{equation}
	\label{MO_model}
	\min_{X} \ (f_0, f_1, \dots, f_{|P|}) ,
\end{equation}
where $X$ is the set of technically feasible fluencies.

This model applies for any method of fluence map optimization, capable to handle functions (\ref{f0}) and (\ref{fp1}).

The classical approach to use model (\ref{MO_model}) for RT planning, with no reference to the gEUD, would be to produce a set of fluence maps $x$ for which $f_0(x),\dots,f_{|P|}(x)$ are mutually nondominated or, with the use of exact optimization methods, even nondominated in $X$. That can be achieved with a form of scalarization of the multi-objective problem (\ref{MO_model}) \citep{EHR10,Kaliszewski_et_all_2016}.

\subsubsection{Combining the multi-objective optimization model and the gEUD-based optimization}
\label{combination}

We propose the following iterative scheme for automated tuning of gEUD parameters, i.e. elements of $\phi$, $\phi \in \Phi$, where $\Phi$ is the set of admissible parameters $\phi$. In each iteration:

\begin{enumerate}
	\item Each collection of gEUD parameters $\phi$ is evaluated with respect to values of $|P + 1|$ criteria $f_0(x^*(\phi)), \dots , f_{|P|}(x^*(\phi))$, where $x^*(\phi) = \argmax_X F(x, \phi)$ is derived by an exact optimizer.
	\item A multi-objective evolutionary optimization algorithm searches set $\Phi$ for collections of parameters $\phi, \ \phi \in \Phi$ , such that tuples $f_0(x^*(\phi)), \dots ,$ $ f_{|P|}(x^*(\phi))$ are mutually nondominated and they dominate analogous tuples derived in the previous iteration.
\end{enumerate} 
The process stops when the differences between tuples in successive iterations become insignificant. The final collection of tuples ($f_0(x^*(\phi)), \dots , f_{|P|}(x^*(\phi))$) represents an approximation of the Pareto efficient set of RT plans $x^*(\phi)$ (the set of all nondominated RT plans at the final iteration).

A distinctive feature of this scheme is that it hybridizes exact and heuristic optimization methods \citep{JOURDAN2009}.

\section{An IMRT planning system with automated gEUD parameter tuning: PersEUD}
\label{sec:perseud}

\begin{figure}[ht]
    	\centering
    	\includegraphics[width=1\textwidth]{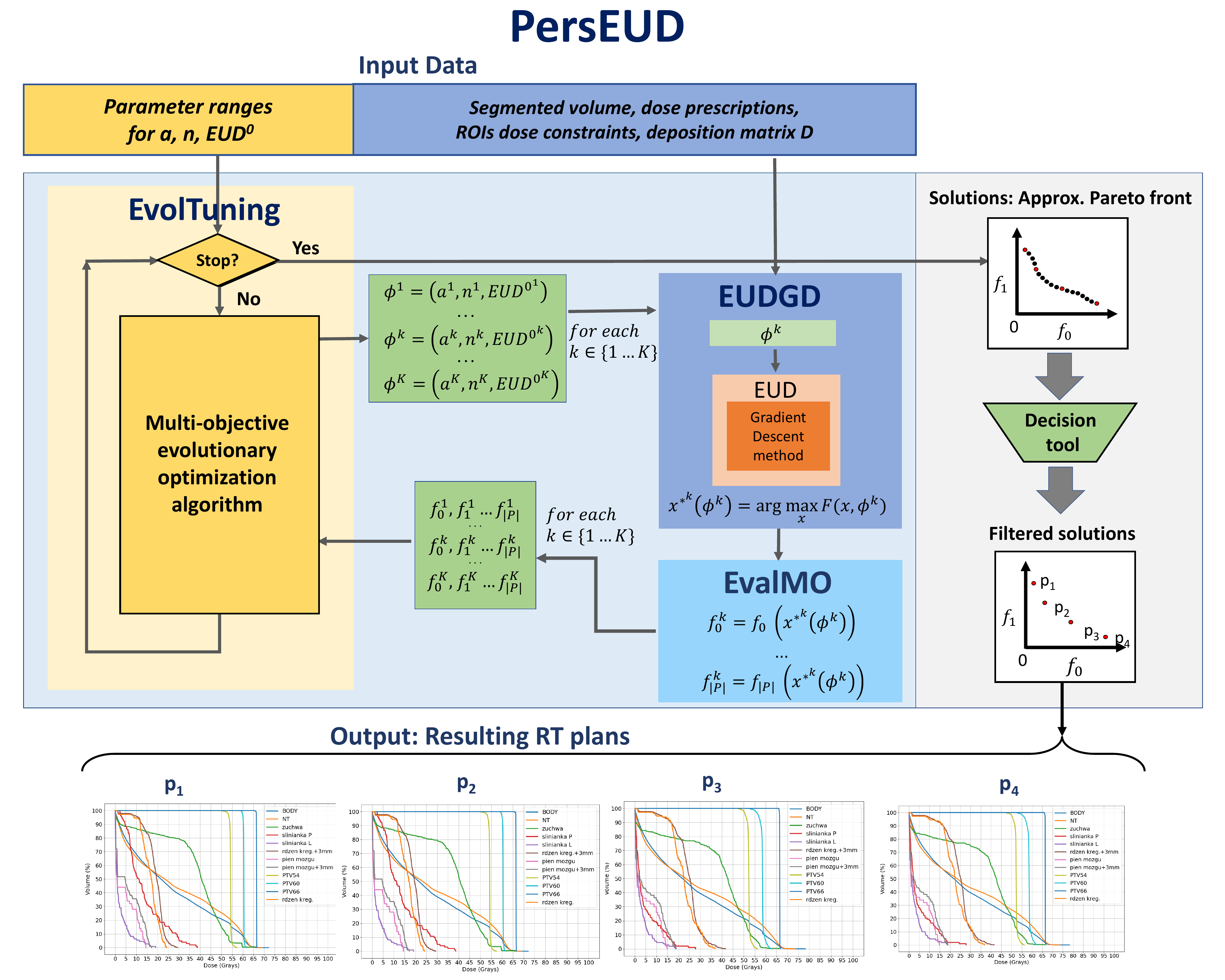}
    	\caption{PersEUD flow diagram.}
    	\label{fig:conceptualModel}
\end{figure}

We have developed and implemented a system, named PersEUD, to deliver IMRT RT plans based on the gEUD, automated parameter $\phi$ tuning, and the optimization hybrid scheme described in Section \ref{combination}. It is composed of four modules:
\begin{enumerate}
\item \textbf{EUDGD}, to deliver the optimal solutions $x^*(\phi)$ maximizing function (\ref{EUDmodel}) for a given collection of parameters $\phi$, 

\item \textbf{EvalMO}, to compute values of functions $f_0, f_1, \dots , f_{|P|}$ for $x^*(\phi)$,

\item \textbf{EvolTuning}, a multi-objective evolutionary algorithm, to provide collections of parameters $\phi$, nondominated in terms of functions $f_0, f_1, \dots , f_{|P|}$,

\item \textbf{Decision Tool}, to reduce the number of collections of parameters $\phi$ presented to RT planners.
\end{enumerate}

The flow diagram of PersEUD is shown in Figure 1. One run of the EvolTuning module produces $K$ collections of parameters $\phi$, mutually nondominated in terms of functions $f_0, f_1, \dots, f_{|P|}$. For each $\phi^k, \ k \in K$, the module EUDGD computes $x^*(\phi^k)$ that maximizes function (\ref{EUDmodel}) and next the module EvalMO computes values of functions $f_0, f_1, \dots, f_{|P|}$. These values are used by the EvolTuning module to identify mutually nondominated parameters $\phi$ directing in turn the search for new $K$ collections of $\phi$. 

It is worth mentioning the versatility of PersEUD, as it can easily accommodate alternative criteria, because criteria $f_0, f_1, \dots, f_{|P|}$ are only involved in EvalMO module. The next four subsections provide more detailed descriptions of the modules of PersEUD.

\subsection{The EUDGD module}

To maximize function (\ref{EUDmodel}), the Gradient Descent method has been custom-coded and the resulting algorithm (GD) is the backbone of the EUDGD module.

The EUDGD module receives as inputs for each patient case: the patient data, the deposition matrix, the beamlet geometry, the ROI segmentation, parameters $\phi$ and the number of steps. It returns the optimal fluence $x^*(\phi)$.
To generate feasible fluence maps, the EUDGD module needs additional configuration parameters. In the current software release these parameters are set to default values, but they can be customized. The default GD parameter values reported in this work are: $20000$ descent steps, $2e^{-7}$ step size, and $0.3$ as the maximum beamlet intensity limit. Lastly, a 3x3 smoothing kernel has been applied to all beamlets after every descent step to facilitate the delivery of fluencies defined by RT plans by radiation equipment.

The number of descent steps required by EUDGD to provide acceptable fluence maps and the size of the auxiliary data structures result in high computational and memory demands. To maximize function (\ref{EUDmodel}) in reasonable times, we have applied HPC software and hardware techniques, as described in \cite{EUDPPAM_22}.

\subsection{The EvalMO module}

For fluence maps $x^*(\phi^k)$ delivered by the EUDGD module, the EvalMO module computes the values of functions $f_0(x^*(\phi^k)), \dots, f_{|P|}(x^*(\phi^k))$ and sends them to EvolTuning module.

\subsection{The EvolTuning module}

The EvolTuning consists of a multi-objective evolutionary optimization algorithm that derives parameters $\phi$. It is initialized with collections of $n_s, a_s, EUD^0_s$, $s \in T \cup R$. 
In order to generate clinically acceptable plans and reduce the search space of the evolutionary process, feasibility ranges for each kind of parameter are to be set accordingly.

For the current implementation, the MOEA/D algorithm has been selected and seeded with parameter values recommended in \cite{zhang2007moea}, with the exception of the number of generations (epochs) that is set to $50$ and the cardinality of populations set to $150$. 

EvolTuning uses the evaluations received from EvalMO to generate new collections of (potentially improved) parameters $\phi$. The process stops when the limit of generations is reached. As the result, $150$ mutually nondominated collections $\phi$ of parameters of function (\ref{EUDmodel}), each collection yielding radiation plan in the form of fluence map $x^*(x)$, 150 plans in all, constitute an approximation of the set of Pareto efficient collections $\phi$ (in terms of functions $f_0, \dots, f_{|P|}$).

\subsection{Decision Tool}
\label{ssec:DecisionTool}

To avoid overwhelming the RT planner with so many plans, the module Decision Tool reduces the set of mutually nondominated collections $\Phi$. It performs this task as follows:

First, plans that minimize objective functions $f_0, f_1, \dots , f_{|P|}$ separately are selected, and they define the extreme points of the respective Pareto front. 
Next, a limited number of plans, fairly distributed between those extreme points, are added to form the reduced approximation of the Pareto front. Observe that even if $f_0$ for some plan takes a positive but small value (meaning that there are some constraint violations), it can be still a clinically viable option as long as it offers plausible values of $f_1, \dots , f_{|P|}$ at the price of a small constraint violation.

\section{Results}
\label{Sec_Results}
As the proof of concept of our proposed methodology, we have prepared three RT plans for real patient data specified in Section \ref{ss:patientdata}. 

The EUDGD module uses a radiation dose deposition model, yielding the deposition matrix $D$,  developed by researchers at the Faculty of Electrical Engineering of the Warsaw University of Technology \citep{wut2015}. The final plan evaluations have been carried out using the commercially available treatment planning system Eclipse (v 15.6, Varian Medical Company). By this, the dose deposition model, as well as the proposed methodology, were tested in a clinical regime.

\subsection{Implementation and experimentation platforms}

The gEUD-based RT planning system, described in the previous section, consists of multiple modules connected by a messaging system based on the ZeroMQ library \citep{ZeroMQ}.

The EvolTuning module is implemented in Java and it takes advantage of the well-known jMetal framework \citep{jmetal2011}. This framework provides well-crafted implementations of multiple algorithms, allowing the user to select the best one for a given problem. Although in this experimentation we have used the MOEA/D algorithm, this module system allows us to swap the evolutionary algorithm with minimal changes.

For the EUDGD and EvalMO modules, we have developed two versions: a CUDA C \citep{CUDA} version for NVIDIA-based GPUs, and an OpenMP C \citep{OpenMP} version for multicore CPUs. These two versions can be used interchangeably or even in parallel to fully exploit all the available resources of the given platform.

Finally, the Decision Tool is a Python-based script that is called after the EvolTuning process terminates. All the figures shown in this section are generated using Python tools implemented in this module.

Our experiments have been run on the High-Performance Cluster of the SAL (Supercomputación--Algoritmos) research group, located at the University of Almeria. Two kinds of computing nodes have been used. The first node contains two AMD EPYC 7302 (32 CPU cores), 512 GB of DDR4 RAM, and two NVIDIA Tesla V100 (32 GB). The second node contains two Intel Xeon E5-2620v3 (12 CPU cores), 64 GB of DDR3 RAM, and two NVIDIA Kepler K80 (12 GB).

\subsection{Patient data and plan evaluation metrics}
\label{ss:patientdata}

For the experimentation, we have used a Head and Neck IMRT real patient case treated with nine radiation beams. In this case, three PTVs with different prescribing radiation dose deposition levels 66 Gy, 60 Gy and 54 Gy (Gray, symbol Gy, is a unit of radiation dose) are defined. The most important is the PTV with the highest prescribed dose (PTV66) because the highest concentration of tumor cells is there. The other two PTVs are treated prophylactically, so the dose distribution homogeneity in them is less critical. In addition, five OARs are singled out: 
the spinal cord +3mm, the brainstem +3mm, the left salivary gland, the right salivary gland, and the mandible. Also, the normal tissue is defined as a region of the patient body outside all OARs and all PTVs. The dose deposition model contains 30265 beamlets interacting with 94647 voxels representing the irradiated part of the patient body. 

Physical restrictions imposed on RT plans by  oncology physicians are converted to the bounds defined in the model presented in Subsection \ref{ssec:multiobjectivemodel}. Table~\ref{t:plan1_results} presents the respective bound values for organs considered in the experiments. Those values are repeated also in Table~\ref{t:plan2_results} and Table~\ref{t:plan3_results}.
The bound values for PTVs are set by the following rules: 

$LB_t = 0.90 \times \textit{prescribed dose for PTV$_t$}$, 

$\overline{LB}_t = 0.98 \times \textit{prescribed dose for PTV$_t$}$, 

$\overline{UB}_t = 1.02 \times \textit{prescribed dose for PTV$_t$}$, 

$UB_t = 1.10 \times \textit{prescribed dose for PTV$_t$}, t \in T$.

Explicit constraints on DVH shapes are not a part of model (\ref{MO_model}) as they are implicitly controlled by gEUD. On the other hand, the DVH shape for PTVs ideally should look like the sequence of three-line segments: horizontal (at 100\% level), vertical (at the prescribed deposited dose value), and horizontal again (at 0\% level) (cf. Figure \ref{fig:study1-dvh}, Figure \ref{fig:study2-dvh}, Figure \ref{fig:study3-dvh}). DVHs are primary tool to verify ex-post the homogeneity of dose distributions in PTVs. As such, they are always visually inspected by RT planers and oncology physicians for a holistic evaluation. A proxy measure of PTV homogeneity is the dose-volume metric Dx$\%$. For a given PTV, Dx$\%$ is the minimal dose deposited in x$\%$ of the most irradiated PTV voxels.
 
In our case, metrics D98$\%$ and D2$\%$ are used, namely plans are to satisfy 

D98$\% \geq$ 95$\% \times \textit{prescribed dose for PTV$_t$}$, and 

D2$\% \leq$ 107$\% \times \textit{prescribed dose for PTV$_t$}$ (i.e., in the latter case: the minimal dose deposited in 2$\%$ of the most irradiated PTV voxels should be less or equal than 107$\% \times EUD^0_t$). 

By this, we have the following levels of Dx$\%$

for PTV66 Gy: D98$\% \geq$ 62.70 Gy, D2$\% \leq$ 70.62 Gy,

for PTV60 Gy: D98$\% \geq$ 57.00 Gy, D2$\% \leq$ 64.20 Gy,

for PTV54 Gy: D98$\% \geq$ 51.30 Gy, D2$\% \leq$ 57.78 Gy.

D2$\%$ and D98$\%$ represent two points on the respective DVHs, thus they give a rough characterization of DVHs.

\subsection{Numerical experiments}
\label{ssec:experiments}
The planning system PersEUD is run three times, each time with a different aim. The aim of the first run is to derive a number (150 by default) of mutually nondominated plans that compromise constraint violations versus the sum of the average of radiation doses deposited in the left and the right salivary glands (a bi-objective optimization problem). In the second run, the aim is to derive plans that compromise constraint violations versus doses deposited in the spinal cord + 3mm (a bi-objective optimization problem), whereas, in the third run, the aim is to derive plans that compromise constraint violations versus doses deposited in the left and the right salivary gland versus doses deposited in the spinal cord +3mm (a tri-criteria optimization problem).

The derived plans are evaluated by medical physicist experts (below: the Experts) who, in their daily work, make treatment plans in a commercial RT planning system in a series of try-and-correct interactions. Once they are satisfied with the result, they submit the plan to the oncology physician responsible for the case for approval. If the plan is disapproved, the process is repeated. In complex cases, two or three cycles may be needed to secure the final physician approval.

In the whole process, they deal with one RT plan at a time. The three runs of PersEUD produce 450 plans, fully automatic. For proof-of-concept analyses, the Experts selected one plan from those produced in each PersEUD run and pre-selected by the Decision Tool module.

\subsection{PersEUD run 1: Compromising constraint violations versus the sum of the average of radiation doses deposited in both salivary glands} 
\label{ssec:case1}

To fulfill the aim of the first PersEUD run, the second objective $f_1$ is the sum of the  average radiation doses deposited in both glands, $|P| = 1$, and $f_0$ is defined by formula (\ref{f0}).

\begin{figure}[htbp]
    \centering
    \includegraphics[width=\textwidth]{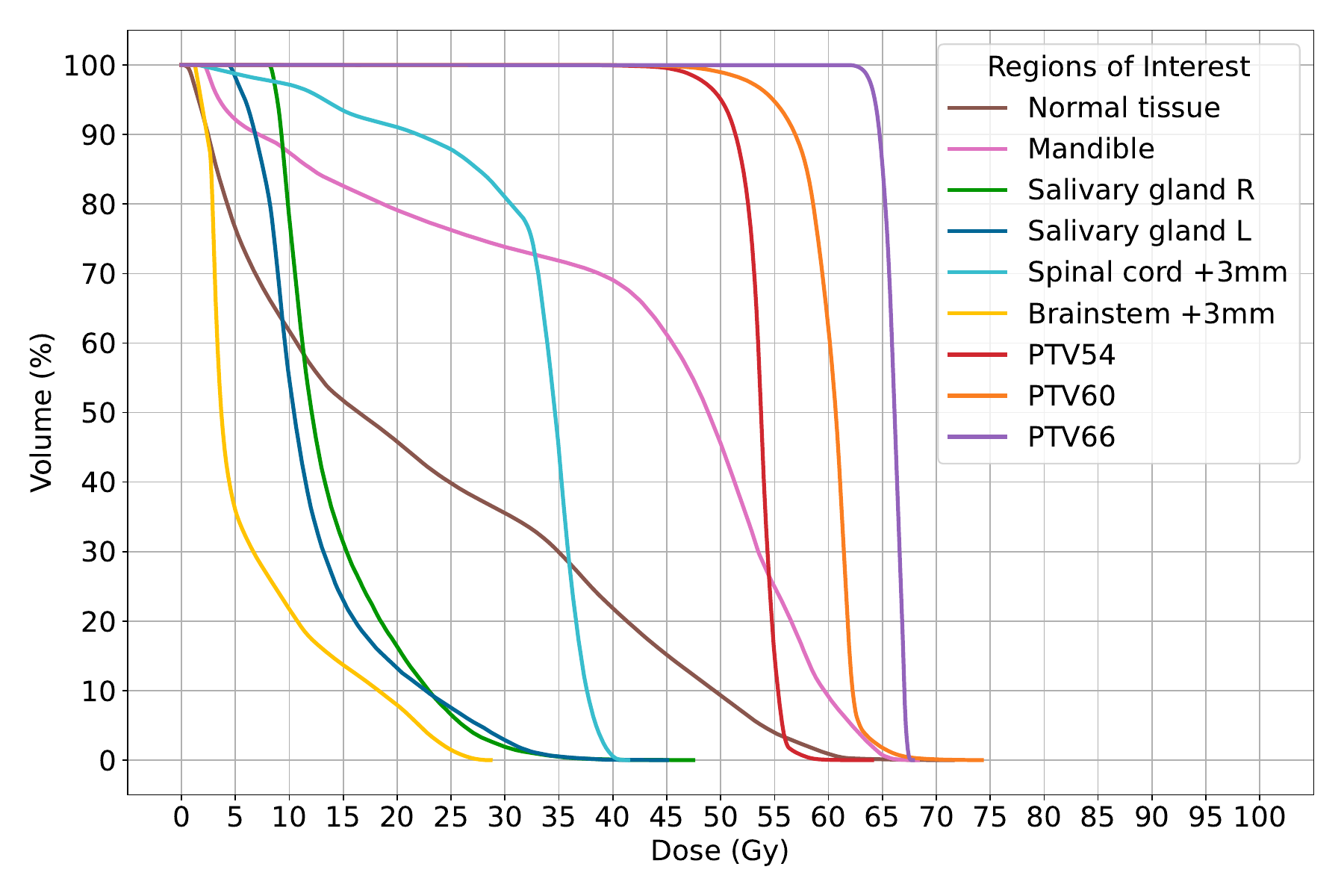}
    \caption{Plan 1. Dose-volume histograms.}
    \label{fig:study1-dvh}
\end{figure}

\begin{figure}[ht]
    \centering
    \begin{subfigure}[b]{0.45\textwidth}
        \centering
        \includegraphics[width=\textwidth]{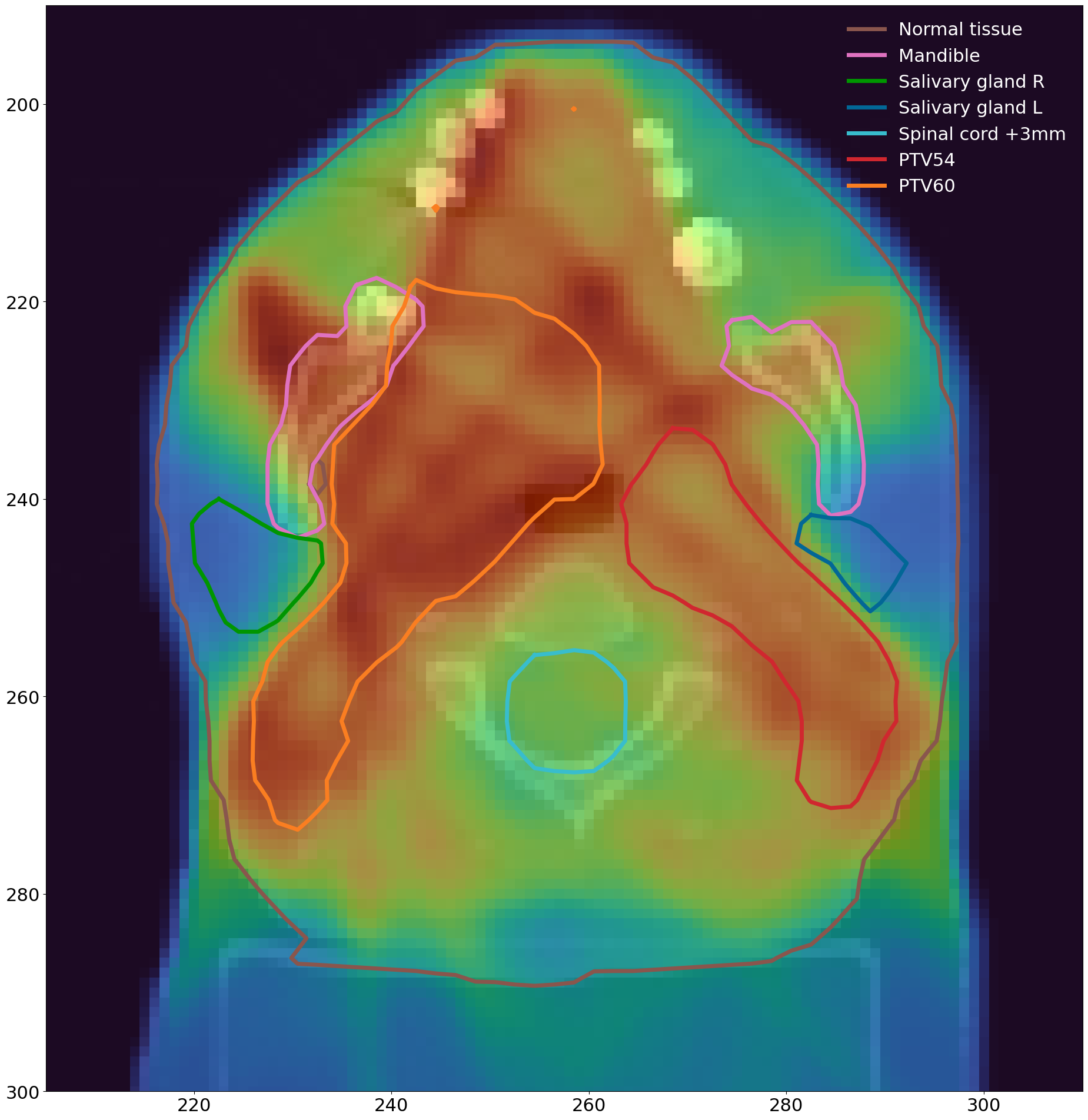}
    \end{subfigure}
    \quad
    \begin{subfigure}[b]{0.45\textwidth}  
        \centering 
        \includegraphics[width=\textwidth]{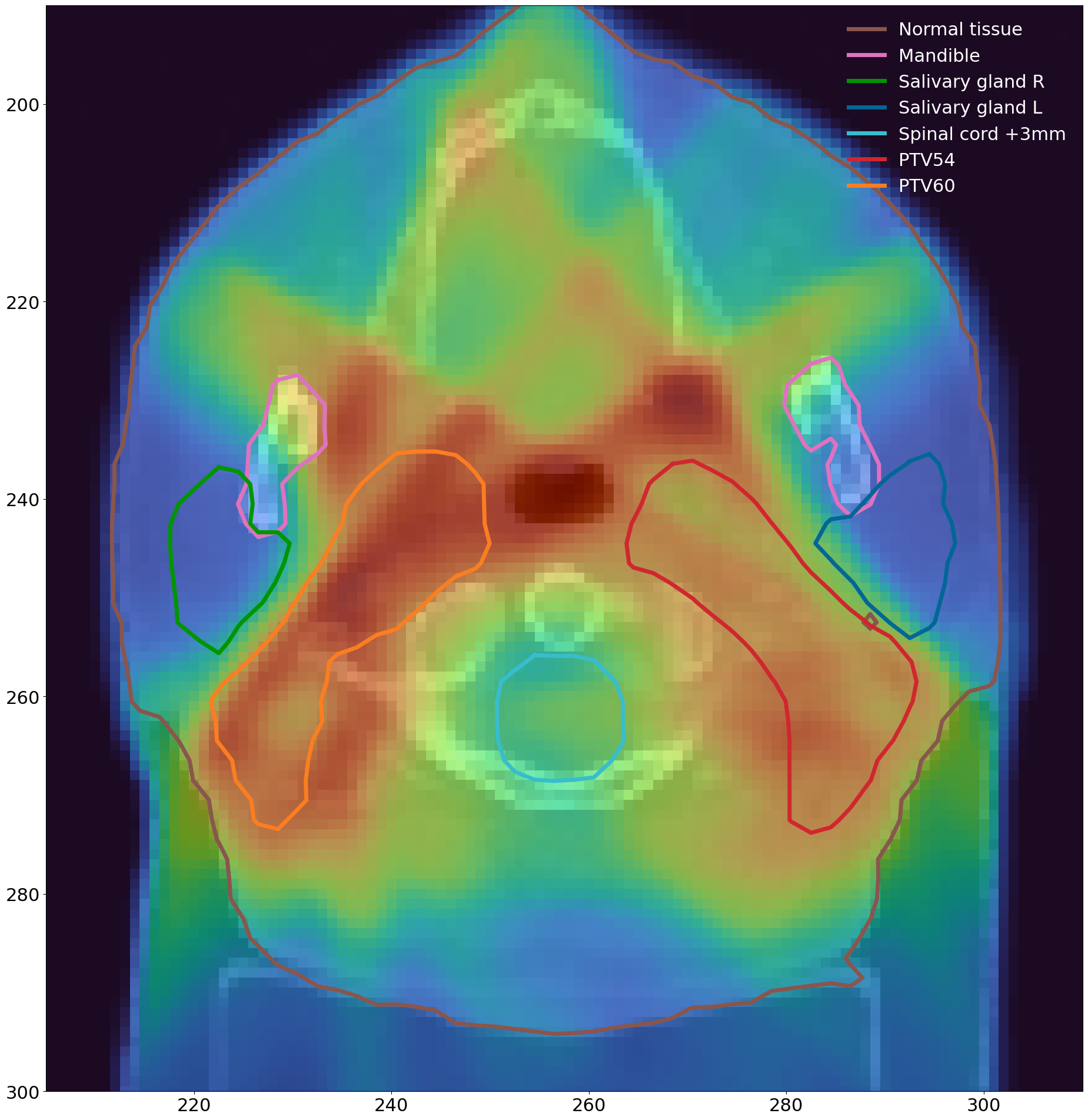}
    \end{subfigure}
    \begin{subfigure}[b]{0.04\textwidth}  
        \centering 
        \includegraphics[width=\textwidth]{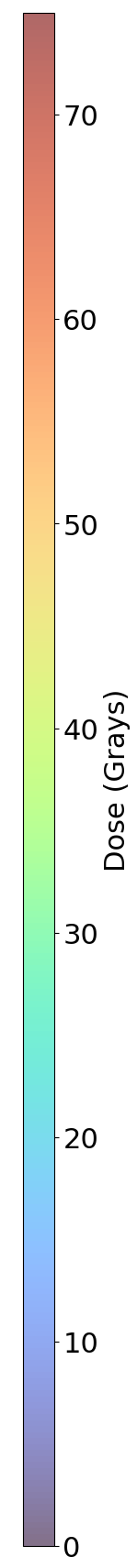}
    \end{subfigure}
    \caption{Plan 1. Dose distributions for $Z = 75$ (left) and $Z=91$ (right). Brown: Normal tissue. Pink: Mandible. Green: Salivary gland R. Blue: Salivary gland L. Cyan: Spinal cord +3mm. Red: PTV54. Orange: PTV60.}
    \label{fig:study1-dose}
\end{figure}

From the resulting RT plans,  Experts select one that does not violate any constraint ($f_0(x) = 0$). This plan is denoted as Plan 1 and it is presented in Table \ref{t:plan1_results}. Table \ref{t:EUDsalivary} presents the parameters of gEUD measures and of the function $F(x, \phi)$ for this plan,
Figure \ref{fig:study1-dvh} presents its dose-volume histograms, and Figure \ref{fig:study1-dose} the dose distributions for this plan on two exemplary cross sections of the patient irradiated part (head and neck).

\begin{table}[htbp]
	\centering
	\caption{Plan 1. Dose bounds, actual doses, and Dx$\%$ metrics. Values that exceed their constraint are marked in bold.}
	\renewcommand{\arraystretch}{1}
	\setlength{\tabcolsep}{5pt}
	\begin{tabular}{lrrrrrrr}
		\toprule
		  & \multicolumn{4}{c}{\textbf{Dose bounds}} & \textbf{Dose} & \multicolumn{2}{c}{\textbf{Dx$\%$}}\\
		\cmidrule{2-5}\cmidrule{7-8} 
		 \textbf{Region of Interest} &  $LB$ &  $\overline {LB}$ &  $\overline{UB}$ &  $UB$ & Act. & Bound & Act. \\
		 \cmidrule{2-8}
		 & \multicolumn{7}{c}{Gy}  \\
		\midrule
		 Normal tissue & - & - & - & 74.25 & 71.52 & - & -\\
		 Mandible & - & - & - & 70.00 & 68.28 & - & -\\
		 Salivary gland R. & - & - & 26.00 & - & 14.26 & - & -\\
		 Salivary gland L. & - & - & 26.00 & - & 12.56 & - & -\\
		 Spinal cord +3mm & - & - & - & 50.00  & 41.43 & - & -\\
		 Brainstem +3mm & - & - & - & 60.00  & 28.67 & - & -\\
		 PTV 54 & 48.60 & 52.92 & 55.08 & 59.40  & 53.47 & - & - \\
		 \hspace{0.2cm}D98$\%$ for PTV 54 & - & - & - & -  & - & 51.30 & {\bf 48.00} \\
		 \hspace{0.2cm}D2$\%$ for PTV 54 & - & - & - & -  & - & 57.78 & 56.22\\
		 PTV 60 & 54.00 & 58.80 & 61.20 & 66.00  & 60.04 & - & -\\
		 \hspace{0.2cm}D98$\%$ for PTV 60 & - & - & - & -  & - & 57.00 & {\bf 52.06}\\
		 \hspace{0.2cm}D2$\%$ for PTV 60  & - & - & - & -  & - & 64.20 &  {\bf 64.70} \\
		 PTV 66 & 59.40 & 64.67 & 67.32 & 72.60  &  66.00 & - & - \\
		 \hspace{0.2cm}D98$\%$ for PTV 66 & - & - & - & -  & - & 62.70 & 63.74\\
	 \hspace{0.2cm}D2$\%$ for PTV 66  & - & - & - & -  &  - & 70.62 & 67.32 \\
		\bottomrule
	\end{tabular}
	\label{t:plan1_results}
\end{table}

\begin{table}[ht]
	\centering
	\caption{Plan 1. Parameters of gEUD ($a_s$) and of $F(x, \phi)$ ($n_s$).  Black: 
			Parameters $\phi$ suggested in the literature. Blue: 
			Parameter search ranges. Green: Optimal parameters derived by the EvolTuning module.}
	\begin{tabular}{lrrr}
		\toprule
		\textbf{Region of Int.} & $\mathbf{EUD^{0}_s}$ & $\mathbf{a_s}$ & $\mathbf{n_s}$\\
		\midrule
		Normal tissue & 74.25 & 40.00 & 5.00\\
		Mandible & 70.00 & 10.00 & 5.00\\
		Salivary gland R.& 
		\textcolor{NavyBlue}{$[0.5, 26]$} $:$ \textcolor{ForestGreen}{$4.37$} & 
		\textcolor{NavyBlue}{$[1, 100]$} $:$ \textcolor{ForestGreen}{$1.01$} & 
		\textcolor{NavyBlue}{$[1, 100]$} $:$ \textcolor{ForestGreen}{$100.00$} \\
		Salivary gland L. & 
		\textcolor{NavyBlue}{$[0.5, 26]$} $:$ \textcolor{ForestGreen}{$3.93$} & 
		\textcolor{NavyBlue}{$[1, 100]$} $:$ \textcolor{ForestGreen}{$1.19$} & 
		\textcolor{NavyBlue}{$[1, 100]$} $:$ \textcolor{ForestGreen}{\phantom{0}$18.04$} \\
		Spinal cord +3mm & 50.00 & 10.00 & 5.00\\
		Brainstem +3mm & 60.00 & 10.00 & 5.00\\
		PTV 54 & 54.00 & 
		\textcolor{NavyBlue}{$[-100, -1]$} $:$ \textcolor{ForestGreen}{$-96.04$} & 
		\textcolor{NavyBlue}{$[1, 100]$} $:$ \textcolor{ForestGreen}{$34.41$} \\
		PTV 60 & 60.00 & 
		\textcolor{NavyBlue}{$[-100, -1]$} $:$ \textcolor{ForestGreen}{$-62.81$} & 
		\textcolor{NavyBlue}{$[1, 100]$} $:$ \textcolor{ForestGreen}{$66.77$} \\
		PTV 66 & 66.00 & 
		\textcolor{NavyBlue}{$[-100, -1]$} $:$ \textcolor{ForestGreen}{$-90.32$} & 
		\textcolor{NavyBlue}{$[1, 100]$} $:$ \textcolor{ForestGreen}{$92.63$} \\
		\bottomrule
	\end{tabular}
	\label{t:EUDsalivary}
\end{table}
   
As seen in Table \ref{t:plan1_results}, this plan provides adequate dose distributions in PTVs, and, as intended, preserves both salivary glands, and it does this surprisingly well. The average doses are 14.26 Gy for the right and 12.56 Gy for the left salivary gland. These values are well below the imposed bound of 26 Gy. All other OARs with no protection priority are also below their acceptable maximal doses. The most notable is the maximum dose in the spinal cord +3mm, equal to 41.43 Gy, while the bound is 50 Gy, and the maximum dose in the brainstem +3mm, equal to 28.67 Gy, while the bound is 60.00 Gy. 

The actual average dose deposited in PTV66 (the most important PTV) is 66 Gy, as prescribed, and D98$\% \geq $ 62.70 Gy and D2$\% \leq 70.62 Gy$.

The actual average dose deposited in PTV60 is 60.04 Gy and that is almost the perfect match as this value is  within $\pm$2$\%$ range from the target value. However, D98$\%$ is significantly unmet (52.06 Gy vs 57.00 Gy expected) due to the protection of the salivary gland. D2$\%$ is also exceeded but by a clinically insignificant value.

The actual average dose deposited in PTV54 is 53.47 Gy and that is also almost the perfect match as this value is within  $\pm$2$\%$ range from the target value. Moreover,  D2$\% \leq 57.78$Gy, but D98$\% \not\geq 51.30$.

As mentioned earlier, in regions treated prophylactically (in the considered case PTV 60 and PTV54) slight violations of D98$\%$ and D2$\%$ related constraints are acceptable.

\FloatBarrier

\subsection{PersEUD run 2: Compromising constraint violations versus doses deposited in the spinal cord +3mm}
\label{ssec:case2}

To fulfill the aim of the second PersEUD run, the second objective $f_1$ is the average radiation dose deposited in the spinal cord +3mm, $|P| = 1$, and $f_0$ is defined by formula (\ref{f0}).

\begin{figure}[ht]
    \centering
    \includegraphics[width=\textwidth]{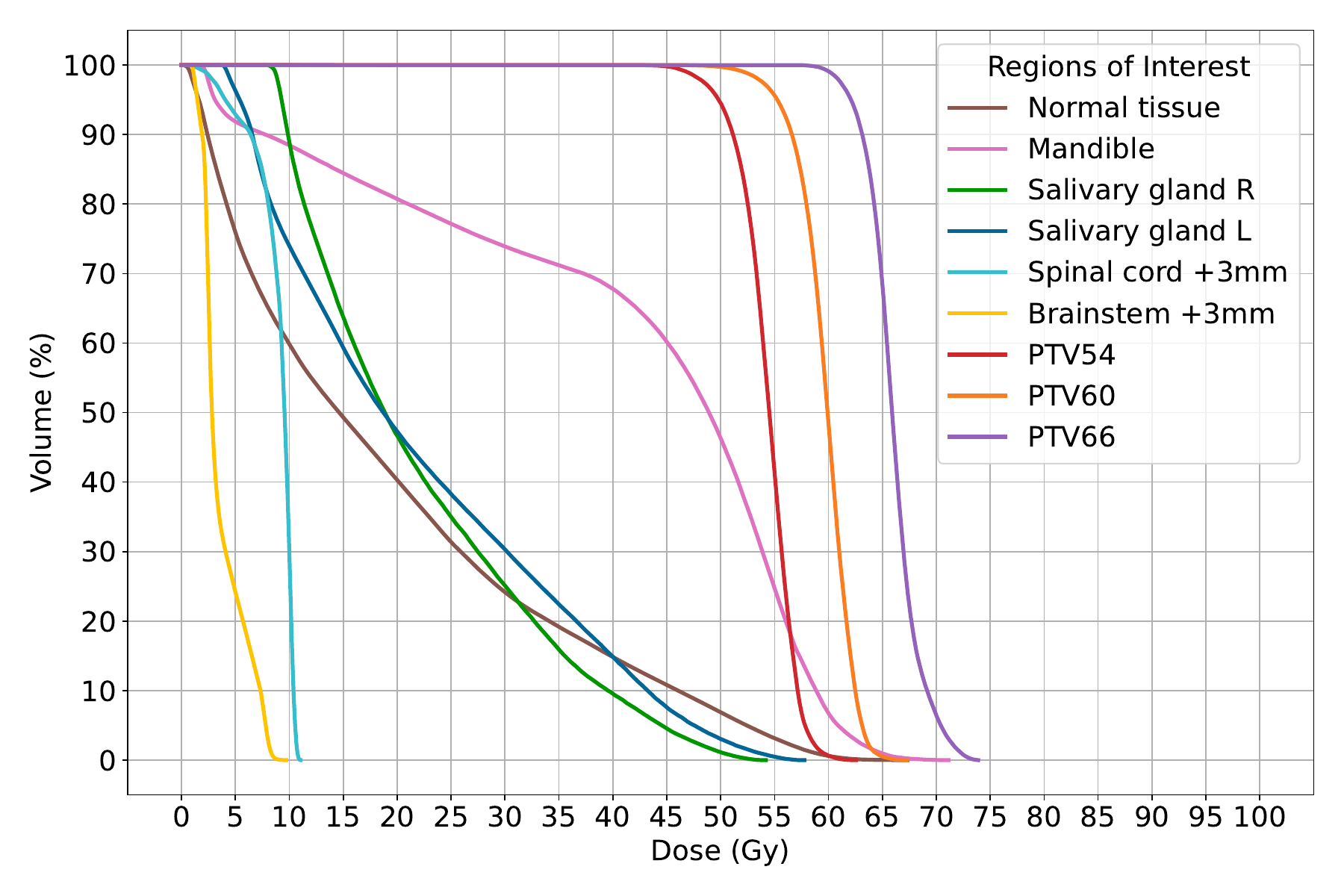}
    \caption{Plan 2. Dose-volume histograms.}    
    \label{fig:study2-dvh}
\end{figure}

\begin{figure}[ht]
    \centering
    \begin{subfigure}[b]{0.45\textwidth}
        \centering
        \includegraphics[width=\textwidth]{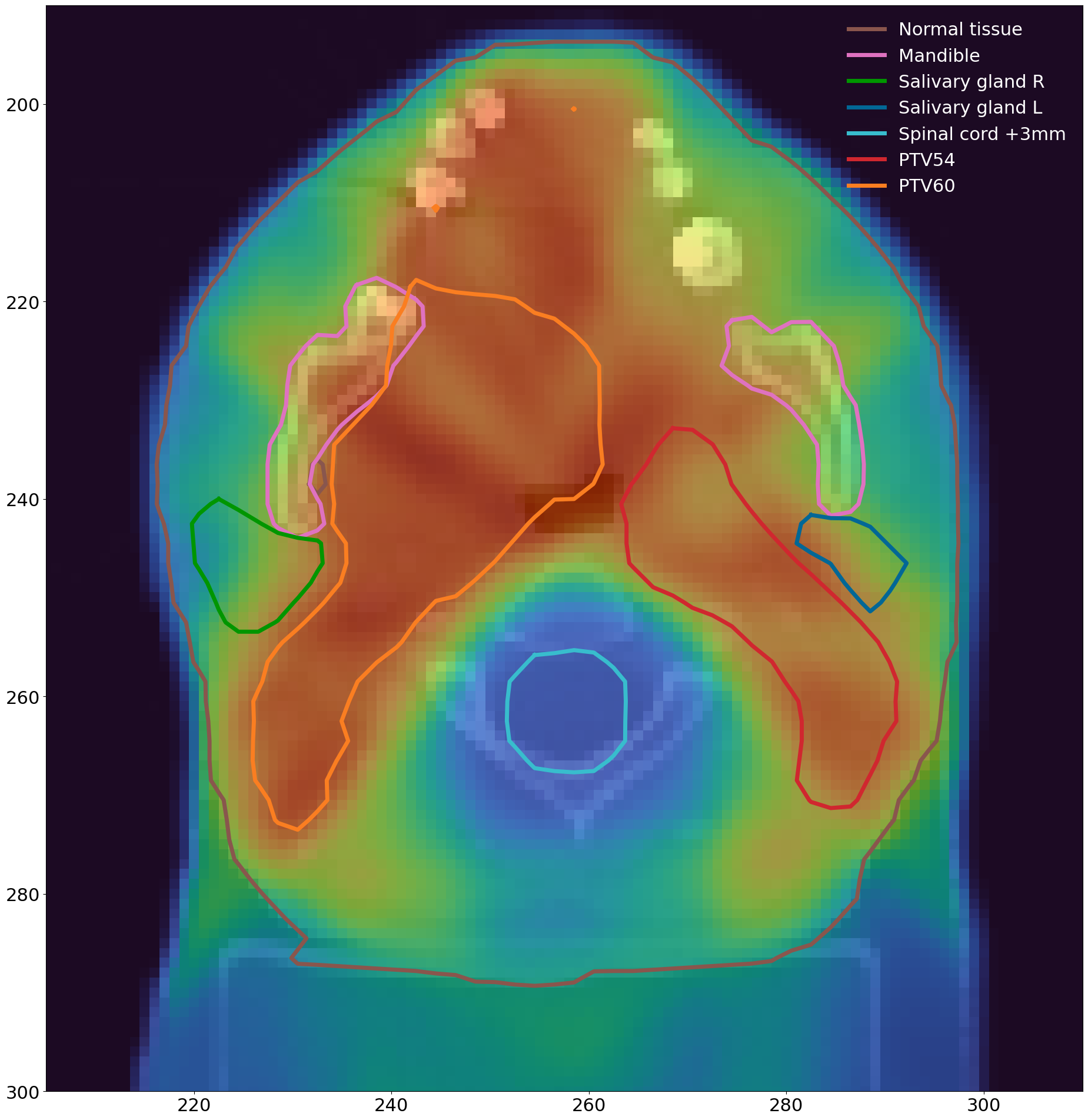}
    \end{subfigure}
    \quad
    \begin{subfigure}[b]{0.45\textwidth}  
        \centering 
        \includegraphics[width=\textwidth]{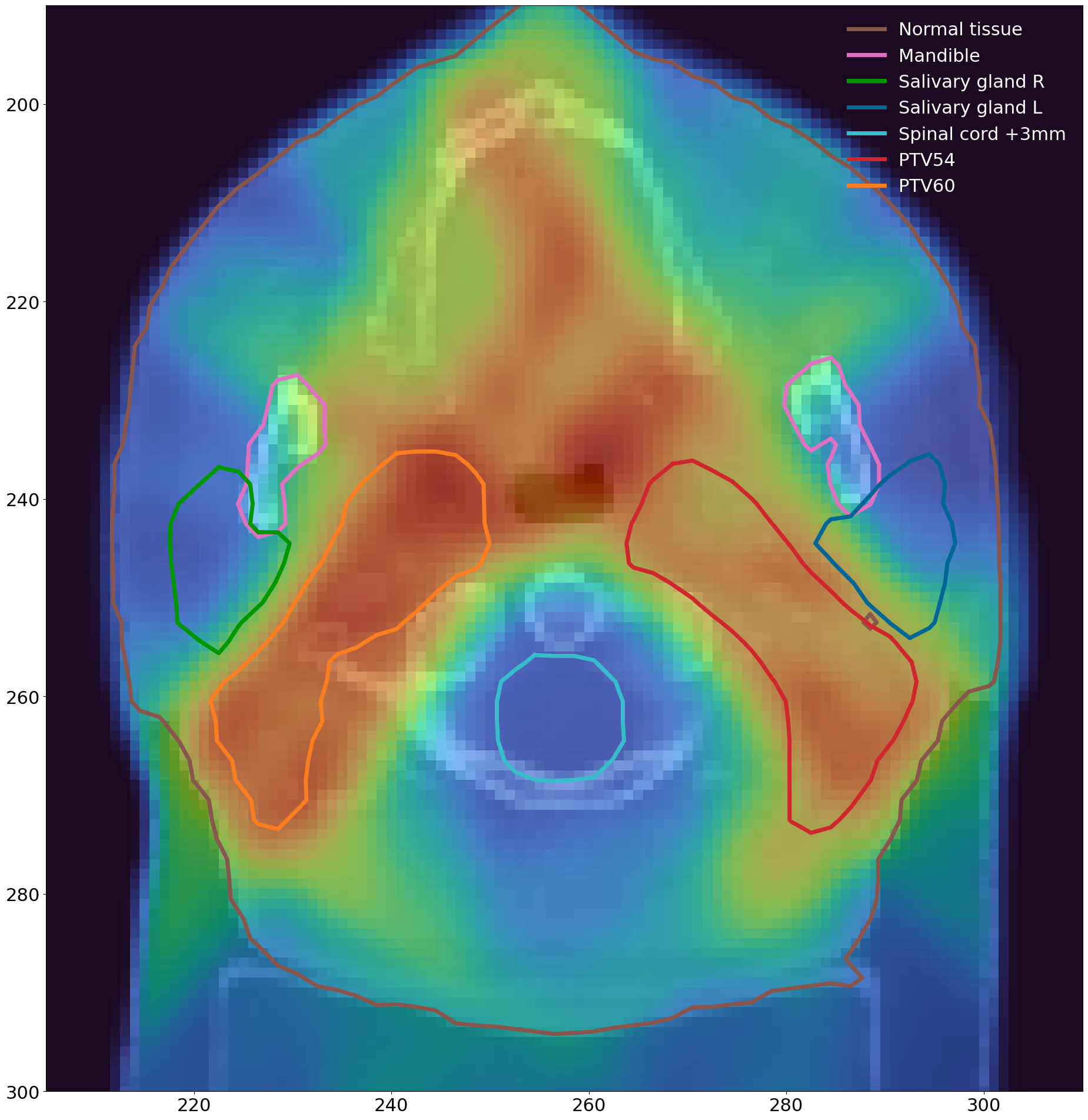}
    \end{subfigure}
    \begin{subfigure}[b]{0.04\textwidth}  
        \centering 
        \includegraphics[width=\textwidth]{colorbar_0_75.png}
    \end{subfigure}
    \caption{Plan 2. Dose distributions for $Z = 75$ (left) and $Z=91$ (right). Brown: Normal tissue. Pink: Mandible. Green: Salivary gland R. Blue: Salivary gland L. Cyan: Spinal cord +3mm. Red: PTV54. Orange: PTV60.}
    \label{fig:study2-dose}
\end{figure}

From the resulting RT plans, the Experts select one that does not violate any constraint ($f_0(x) = 0$). This plan is denoted as Plan 2 and is presented in Table \ref{t:plan2_results}. Table \ref{t:EUDspinal} presents parameters of gEUD and of function $F(x, \phi)$ for this plan, Figure \ref{fig:study2-dvh} presents its dose-volume histograms, and Figure \ref{fig:study2-dose} the dose distributions for this plan on two exemplary cross sections of the patient irradiated part (head and neck).

\begin{table}[htbp]
	\centering
	\caption{Plan 2. Dose bounds, actual doses, and Dx$\%$ metrics. Values that exceed their constraint are marked in bold.}
	\renewcommand{\arraystretch}{1}
	\setlength{\tabcolsep}{5pt}
	\begin{tabular}{lrrrrrrr}
		\toprule
		  & \multicolumn{4}{c}{\textbf{Dose bounds}} & \textbf{Dose} & \multicolumn{2}{c}{\textbf{Dx$\%$}}\\
		\cmidrule{2-5}\cmidrule{7-8} 
		 \textbf{Region of Interest} &  $LB$ &  $\overline {LB}$ &  $\overline{UB}$ &  $UB$ & Act. & Bound & Act. \\
		 \cmidrule{2-8}
		 & \multicolumn{7}{c}{Gy}  \\
		\midrule
		 Normal tissue & - & - & - & 74.25 & 67.30 & - & -\\
		 Mandible & - & - & - & 70.00 & {\bf 71.13} & - & -\\
		 Salivary gland R. & - & - & 26.00 & - & 22.08 & - & -\\
		 Salivary gland L. & - & - & 26.00 & - & 22.19 & - & -\\
		 Spinal cord +3mm & - & - & - & 50.00  & 11.02 & - & -\\
		 Brainstem +3mm & - & - & - & 60.00  & 9.71 & - & -\\
		 PTV 54 & 48.60 & 52.92 & 55.08 & 59.40  & 54.30 & - & - \\
		 \hspace{0.2cm}D98$\%$ for PTV 54 &  &  & - & -  & - & 51.30 & {\bf 48.03} \\
		 \hspace{0.2cm}D2$\%$ for PTV 54 & - & - & - & -  & - & 57.78 & 58.81\\
		 PTV 60 & 54.00 & 58.80 & 61.20 & 66.00  & 59.68 & - & -\\
		 \hspace{0.2cm}D98$\%$ for PTV 60 & - & - & - & -  & - & 57.00 & 58.81\\
		 \hspace{0.2cm}D2$\%$ for PTV 60  & - & - & - & -  & - & 64.20 & 63.80 \\
		 PTV 66 & 59.40 & 64.67 & 67.32 & 72.60  &  66.00 & - & - \\
		 \hspace{0.2cm}D98$\%$ for PTV 66 & - & - & - & -  & - & 62.70 & {\bf 60.87}\\
		 \hspace{0.2cm}D2$\%$ for PTV 66  & - & - & - & -  &  - & 70.62 & {\bf 71.66} \\
		\bottomrule
	\end{tabular}
	\label{t:plan2_results}
\end{table}

\begin{table}[ht]
	\centering
	\caption{Plan 2. Parameters of gEUD ($a_s$) and of $F(x, \phi)$ ($n_s$).  Black: 
			Parameters $\phi$ suggested in the literature. Blue: 
			Parameter search ranges. Green: Optimal parameters derived by the EvolTuning module.}
	
	\begin{tabular}{lrrr}
		\toprule
		\textbf{Region of Int.} & $\mathbf{EUD^{0}_s}$ & $\mathbf{a_s}$ & $\mathbf{n_s}$\\
		\midrule
		Normal tissue & 74.25 & 40.00 & 5.00\\
		Mandible & 70.00 & 10.00 & 5.00\\
		Salivary gland R.& 26.00 & 1.00 & 5.00\\
		Salivary gland L. & 26.00 & 1.00 & 5.00 \\
		Spinal cord +3mm & 
		\textcolor{NavyBlue}{$[0.5, 50]$} $:$ \textcolor{ForestGreen}{$0.50$} & 
		\textcolor{NavyBlue}{$[1, 50]$} $:$ \textcolor{ForestGreen}{$49.30$} & 
		\textcolor{NavyBlue}{$[1, 100]$} $:$ \textcolor{ForestGreen}{$10.25$} \\
		Brainstem +3mm & 60.00 & 10.00 & 5.00\\
		PTV 54 & 54.00 & 
		\textcolor{NavyBlue}{$[-100, -1]$} $:$ \textcolor{ForestGreen}{$-67.50$} & 
		\textcolor{NavyBlue}{$[1, 100]$} $:$ \textcolor{ForestGreen}{$17.76$} \\
		PTV 60 & 60.00 & 
		\textcolor{NavyBlue}{$[-100, -1]$} $:$ \textcolor{ForestGreen}{$-49.83$} & 
		\textcolor{NavyBlue}{$[1, 100]$} $:$ \textcolor{ForestGreen}{$25.35$} \\
		PTV 66 & 66.00 & 
		\textcolor{NavyBlue}{$[-100, -1]$} $:$ \textcolor{ForestGreen}{\phantom{0}$-7.02$} & 
		\textcolor{NavyBlue}{$[1, 100]$} $:$ \textcolor{ForestGreen}{$14.07$} \\
		\bottomrule
	\end{tabular}
	
	\label{t:EUDspinal}
\end{table}

\FloatBarrier

In Plan 2 the maximum dose deposited in the spinal cord +3mm is 11.02 Gy, noticeably lower than in the Plan 1. Although protection of the salivary glands in this plan is not a priority, the average dose is 22.19 Gy and 22.08 Gy in the left and right salivary gland, respectively, which is far below the dose tolerance for these organs.
This is at the price of a slight violation of the maximal dose allowed in the mandible (70 Gy) where the actual dose deposited is 71.13 Gy.

The actual average dose deposited in PTV66 is 66 Gy, as prescribed.	However, D98$\% \not\geq $ 62.70 Gy (violation by 1.83 Gy), and D2$\% \not\leq$ 70.62 Gy (violation by 1.04 Gy).

The actual average dose deposited in PTV60 is 59.68 Gy and that is almost the perfect match as this value is  within $\pm$2$\%$ range from the target value. Moreover,  D98$\% \geq$ 57.00 Gy and D2$\% \leq $ 64.20 Gy.

The actual average dose deposited in PTV54 is 54.30 Gy and that is also almost the perfect match as this value is within  $\pm$2$\%$ range from the target value. However, D98$\% \not\geq$ 51.30 Gy, but D2$\% \leq $ 57.78 Gy.

As mentioned earlier, in regions treated prophylactically (i.e. PTV 60 and PTV54) slight violations of D98$\%$ and D2$\%$ related constraints are acceptable.

\FloatBarrier
 
\subsection{PersEUD run 3: Compromising constraint violations versus doses deposited in the left and the right salivary gland versus doses deposited in the spinal cord +3mm}	
\label{ssec:case3}

To fulfill the aim of the third PersEUD run, the second objective function $f_1$ is the average of the radiation doses deposited in the left and the right salivary glands, the third objective function the is the maximal radiation dose deposited in the spinal cord +3mm, $|P| = 2$, and $f_0$ is defined by formula (\ref{f0}).

\begin{table}[ht]
	\centering
	\caption{Plan 3. Dose bounds, actual doses, and Dx$\%$ metrics. Values that exceed their constraint are marked in bold}
	\renewcommand{\arraystretch}{1}
	\setlength{\tabcolsep}{5pt}
	\begin{tabular}{lrrrrrrr}
		\toprule
		  & \multicolumn{4}{c}{\textbf{Dose bounds}} & \textbf{Dose} & \multicolumn{2}{c}{\textbf{Dx$\%$}}\\
		\cmidrule{2-5}\cmidrule{7-8} 
	 \textbf{Region of Interest} &  $LB$ &  $\overline {LB}$ &  $\overline{UB}$ &  $UB$ & Act. & Bound & Act. \\
	 \cmidrule{2-8}
	 & \multicolumn{7}{c}{Gy}  \\
		\midrule
		 Normal tissue & - & - & - & 74.25 & 73.14 & - & -\\
		 Mandible & - & - & - & 70.00 & 69.83 & - & -\\
		 Salivary gland R. & - & - & 26.00 & - & 13.94 & - & -\\
		 Salivary gland L. & - & - & 26.00 & - & 13.14 & - & -\\
		 Spinal cord +3mm & - & - & - & 50.00  & 24.06 & - & -\\
		 Brainstem +3mm & - & - & - & 60.00  & 21.78 & - & -\\
		 PTV 54 & 48.60 & 52.92 & 55.08 & 59.40  & 55.34 & - & - \\
		 \hspace{0.2cm}D98$\%$ for PTV 54 & - & - & - & -  & - & 51.30 & {\bf 49.64} \\
		 \hspace{0.2cm}D2$\%$ for PTV 54 & - & - & - & -  & - & 57.78 & {\bf 59.13}\\
		 PTV 60 & 54.00 & 58.80 & 61.20 & 66.00  & 60.27 & - & -\\
		 \hspace{0.2cm}D98$\%$ for PTV 60 & - & - & - & -  & - & 57.00 & {\bf 52.30}\\
		 \hspace{0.2cm}D2$\%$ for PTV 60  & - & - & - & -  & - & 64.20 &  {\bf 65.49} \\
		 PTV 66 & 59.40 & 64.67 & 67.32 & 72.60  &  66.00 & - & - \\
		 \hspace{0.2cm}D98$\%$ for PTV 66 & - & - & - & -  & - & 62.70 & 63.15 \\
		 \hspace{0.2cm}D2$\%$ for PTV 66  & - & - & - & -  &  - & 70.62 & 67.63 \\
		\bottomrule
	\end{tabular}
	\label{t:plan3_results}
\end{table}

\begin{table}[ht]
	\centering
	\caption{Plan 3. Parameters of gEUD ($a_s$) and of $F(x, \phi)$ ($n_s$).  Black: 
			Parameters $\phi$ suggested in the literature. Blue: 
			Parameter search ranges. Green: Optimal parameters derived by the EvolTuning module.}
	\begin{tabular}{lrrr}
		\toprule
		\textbf{Region of Int.} & $\mathbf{EUD^{0}_s}$ & $\mathbf{a_s}$ & $\mathbf{n_s}$\\
		\midrule
		Normal tissue & 74.25 & 40.00 & 5.00\\
		Mandible & 70.00 & 10.00 & 5.00\\
		Salivary gland R.& 
		\textcolor{NavyBlue}{$[0.5, 26]$} $:$ \textcolor{ForestGreen}{$4.76$} & 
		\textcolor{NavyBlue}{$[1, 100]$} $:$ \textcolor{ForestGreen}{$1.01$} & 
		\textcolor{NavyBlue}{$[1, 100]$} $:$ \textcolor{ForestGreen}{$18.25$} \\
		Salivary gland L. & 
		\textcolor{NavyBlue}{$[0.5, 26]$} $:$ \textcolor{ForestGreen}{$3.79$} & 
		\textcolor{NavyBlue}{$[1, 100]$} $:$ \textcolor{ForestGreen}{$1.31$} & 
		\textcolor{NavyBlue}{$[1, 100]$} $:$ \textcolor{ForestGreen}{$11.17$} \\
		Spinal cord +3mm & 
		\textcolor{NavyBlue}{$[0.5, 50]$} $:$ \textcolor{ForestGreen}{$1.80$} & 
		\textcolor{NavyBlue}{$[1, 50]$} $:$ \textcolor{ForestGreen}{$1.33$} & 
		\textcolor{NavyBlue}{$[1, 100]$} $:$ \textcolor{ForestGreen}{$12.79$} \\
		Brainstem +3mm & 60.00 & 10.00 & 5.00\\
		PTV 54 & 54.00 & 
		\textcolor{NavyBlue}{$[-100, -1]$} $:$ \textcolor{ForestGreen}{$-65.55$} & 
		\textcolor{NavyBlue}{$[1, 100]$} $:$ \textcolor{ForestGreen}{$54.11$} \\
		PTV 60 & 60.00 & 
		\textcolor{NavyBlue}{$[-100, -1]$} $:$ \textcolor{ForestGreen}{$-87.12$} & 
		\textcolor{NavyBlue}{$[1, 100]$} $:$ \textcolor{ForestGreen}{$57.66$} \\
		PTV 66 & 66.00 & 
		\textcolor{NavyBlue}{$[-100, -1]$} $:$ \textcolor{ForestGreen}{$-33.27$} & 
		\textcolor{NavyBlue}{$[1, 100]$} $:$ \textcolor{ForestGreen}{$18.62$} \\
		\bottomrule
	\end{tabular}
	\label{t:Case3}
\end{table}  

\FloatBarrier

\begin{figure}[ht]
    \centering
    \includegraphics[width=\textwidth]{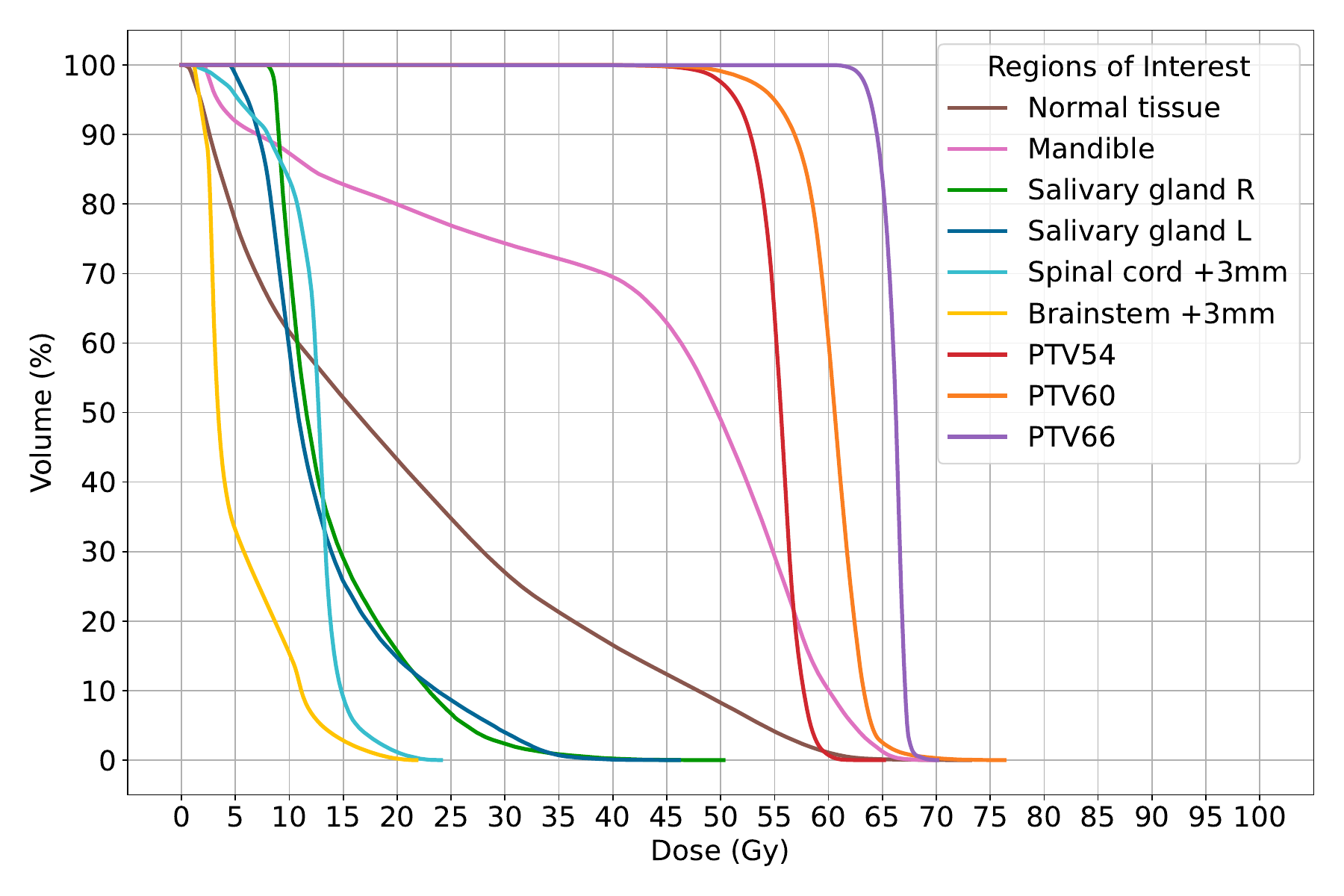}
    \caption{Plan 3. Dose-volume histograms.}    
    \label{fig:study3-dvh}
\end{figure} 

\begin{figure}[ht]
    \centering
    \begin{subfigure}[b]{0.45\textwidth}
        \centering
        \includegraphics[width=\textwidth]{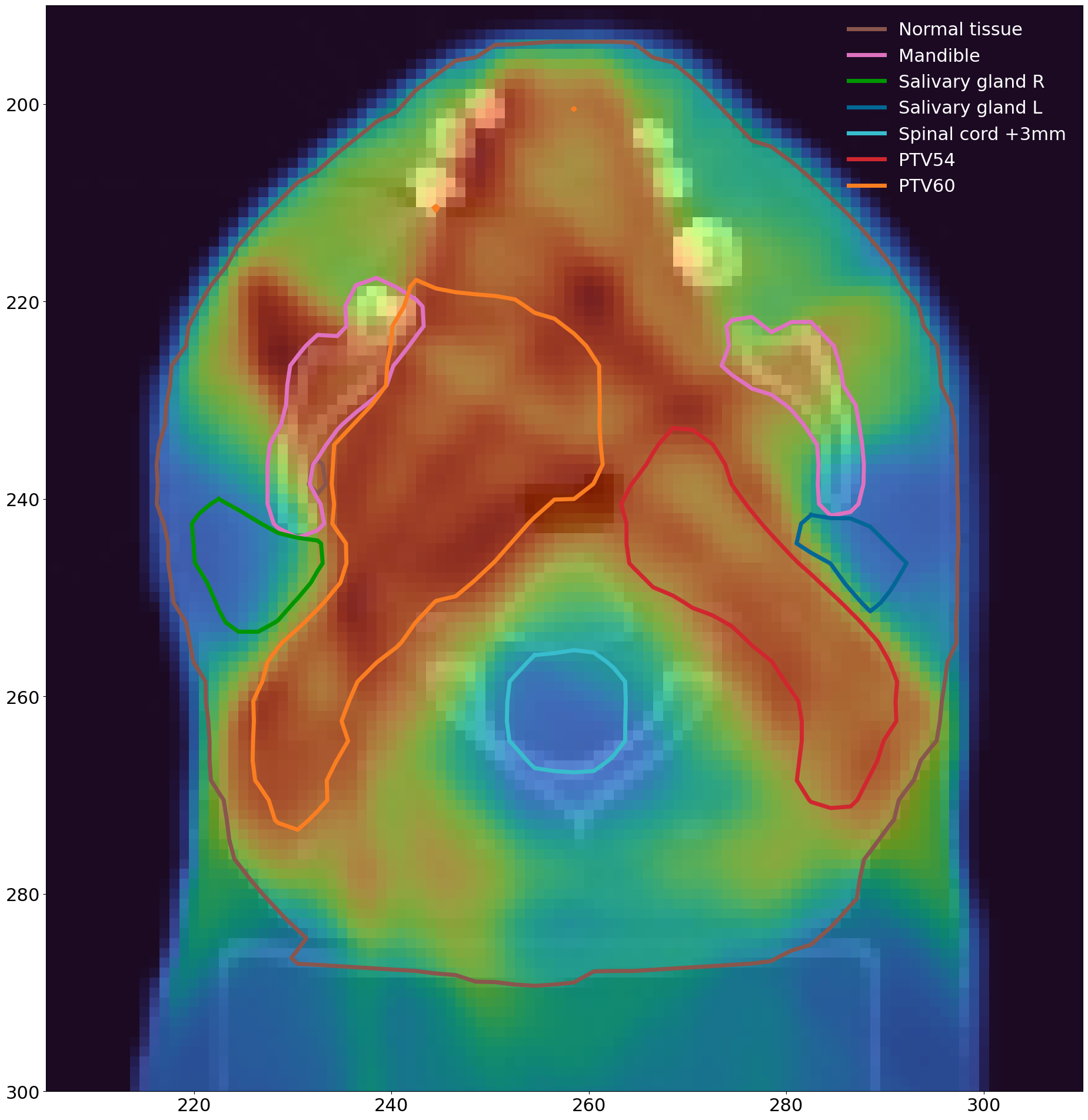}
    \end{subfigure}
    \quad
    \begin{subfigure}[b]{0.45\textwidth}  
        \centering 
        \includegraphics[width=\textwidth]{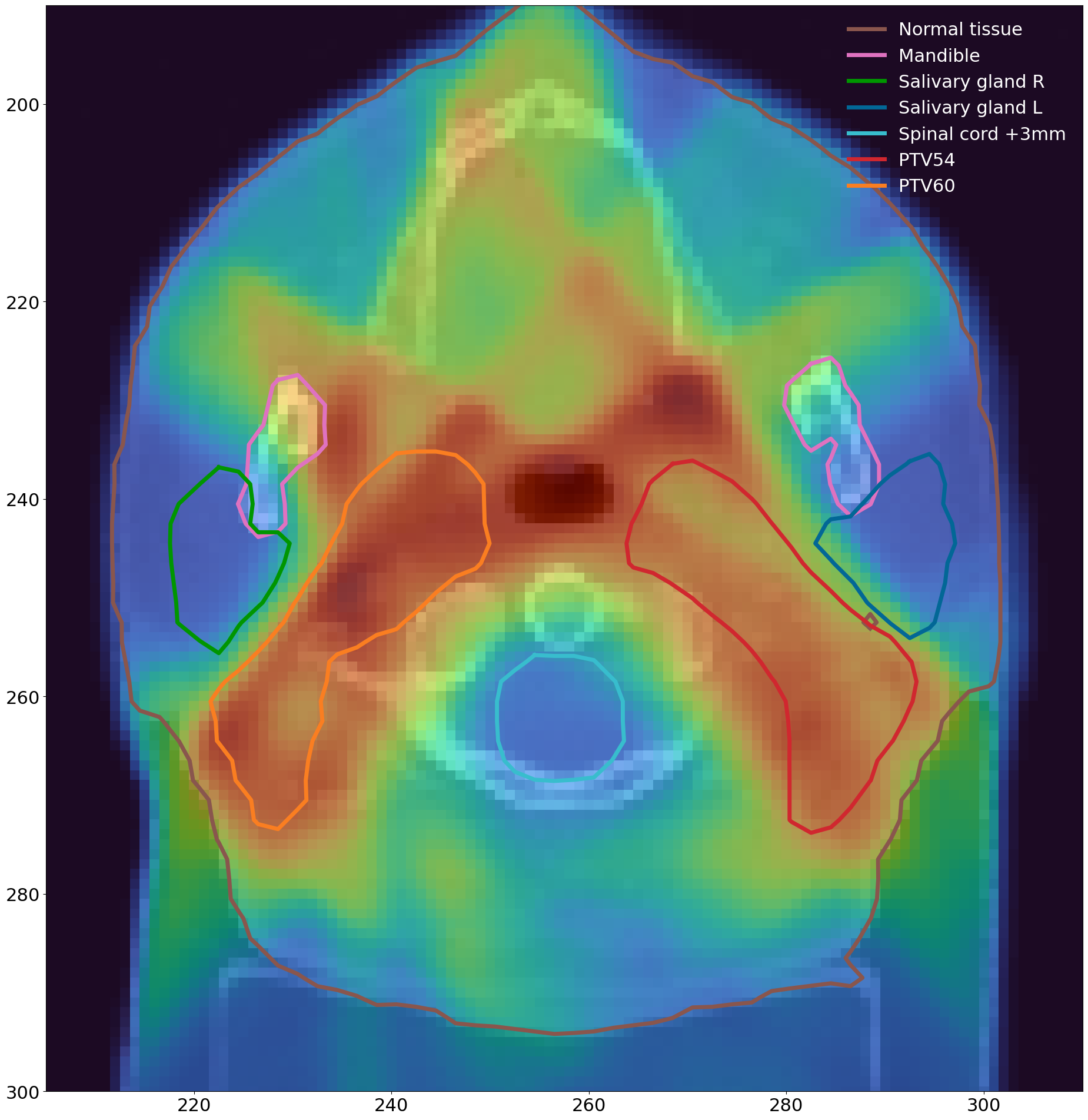}
    \end{subfigure}
    \begin{subfigure}[b]{0.04\textwidth}  
        \centering 
        \includegraphics[width=\textwidth]{colorbar_0_75.png}
    \end{subfigure}
    \caption{Plan 3. Dose distributions for $Z = 75$ (left) and $Z=91$ (right). Brown: Normal tissue. Pink: Mandible. Green: Salivary gland R. Blue: Salivary gland L. Cyan: Spinal cord +3mm. Red: PTV54. Orange: PTV60.}
    \label{fig:study3-dose}
\end{figure}

From the resulting RT plans, the Experts select one that does not violate any constraint ($f(x) = 0$). This plan is denoted as Plan 3 and is presented in Table \ref{t:plan3_results}. Table \ref{t:Case3} presents parameters of gEUD and of function $F(x, \phi)$ for this plan, Figure \ref{fig:study3-dvh} presents its dose-volume histograms, and Figure \ref{fig:study3-dose} the dose distributions for this plan on two exemplary cross sections of the 3D model of the patient irradiated part (head and neck).

In Plan 3, the average dose deposited in the right and the left salivary gland is 13.94 Gy and 13.14 Gy, respectively. That is much lower than the maximal admissible value of 26 Gy, while the maximal dose deposited in the spinal cord +3mm is 24.06 Gy (the maximal acceptable value is 50 Gy). Although Plan 1 and Plan 2 offer slightly lower doses in respective priority protected OARs (the salivary glands in the first case, and the spinal cord +3mm in the second case), Plan 3 well balances the doses delivered to both salivary glands and to the spinal cord +3mm.

The actual average dose deposited in PTV66 is 66 Gy, as prescribed, and D98$\% \geq $ 62.70 Gy, and D2$\% \leq $ 70.62 Gy. The actual average dose deposited in PTV60 is 60.2768 Gy, within $\pm$2$\%$ range from the target value. Moreover,  D98$\% \not\geq $ 57.00 Gy and D2$\% \leq $ 64.20 Gy. The actual average dose deposited in PTV54 is 55.34 Gy and that is also almost the perfect match as this value is within $\pm$2.5$\%$ range from the target value. However, D98$\% \not\geq$ 51.30 Gy, but D2$\% \leq $ 57.78 Gy.

As mentioned earlier, in regions treated prophylactically (i.e. PTV 60 and PTV54) slight violations of D98$\%$ and D2$\%$ related constraints are acceptable.

\subsection{Analysis of PersEUD plans in a clinical setting}
\label{Sec:Analysis}
Table \ref{Cases} provides a comparison of the three plans. At a glance, the Experts disapproved Plan 2 due to its underperformance in PTV66, the most important PTV in the considered case. The explanation for the disapproval is that such an underperformance cannot be outweighed by the excellent performance in the spinal cord +3mm and the brainstem +3mm. From the remaining two, the Experts prefer Plan 3 because it offers low radiation doses deposited in both salivary glands, the spinal cord +3mm, and brainstem +3mm.

\begin{table}[htbp!]
\centering
\caption{A comparison table for the three plans. 
	Values that exceed their corresponding constraint are marked in bold.}
\begin{tabular}{lllrrr}
            \toprule
       \textbf{ROI} & \textbf{Bound} & \textbf{Bound} & \textbf{Plan 1} & \textbf{Plan 2} & \textbf{Plan 3} \\ \cmidrule{3-6}
           & {\bf type}  & \multicolumn{4}{c}{\textbf{Gy}}  \\
            \midrule
       Normal tissue & 
            Average &
            74.25 &
            71.52 & 
            67.30 &
            73.14 \\
       Mandible &
            Maximum &
            70.00 &
            68.28 &
            \bf{71.13} &
            69.83 \\
       Salivary gland R.& 
                Average &
                26.00 &
                14.26 & 
                22.08 &
                13.94\\
       Salivary gland L. & 
                Average   & 
                26.00 &
                12.56 & 
                22.19 &
                13.14\\
       Spinal cord +3mm & 
            Maximum &
            50.00 &
            41.43 &
            11.02&
            24.06 \\
       Brainstem +3mm & 
            Maximum &
            60.00 &
            28.67 &
            9.71 &
            21.79 \\
        PTV 54 & 
             Average  &54.00 
             & 53.47 & 54.30 &55.34 \\
           & $D98\%$ &51.30 &\bf{48.00} & \bf{48.03} & \bf{49.64} \\
           & $D2\%$ & 57.78 & 56.22 & \bf{58.81} & \bf{59.13} \\
        PTV 60 & 
              Average  & 60.00 
              & 60.04 & 59.68 & 60.27 \\
             & $D98\%$ & 57.00 & \bf{52.06} & 58.81 & \bf{52.30} \\
             & $D2\%$ & 64.20 & \bf{64.76} & 63.80 & \bf{65.49} \\
         PTV 66 & 
            Average  & 66.00 & 66.00 & 66.00 & 66.00 \\
         & $D98\%$ & 62.70 & 63.74 & \bf{60.87} & 63.15 \\
         & $D2\%$ & 70.62 & 67.32 & \bf{71.66} & 67.63 \\
            \bottomrule
        \end{tabular}
\label{Cases}
\end{table}

To verify the clinical viability of Plan 3 (below: P3), we compare it with two plans prepared by the Expert in Eclipse\textsuperscript{\tiny TM} Treatment Planning system (below: Eclipse) from Varian Medical Systems, namely:

\begin{itemize}
    \item[{R1}] \textit{Starting plan}: an RT plan generated automatically by Eclipse at the start of the planning process.
    \item[{R2}] \textit{Blind-expert plan}: an RT plan prepared by the Expert without seeing P3.
\end{itemize}

\begin{figure}[ht]
    \centering
    \includegraphics[width=\textwidth]{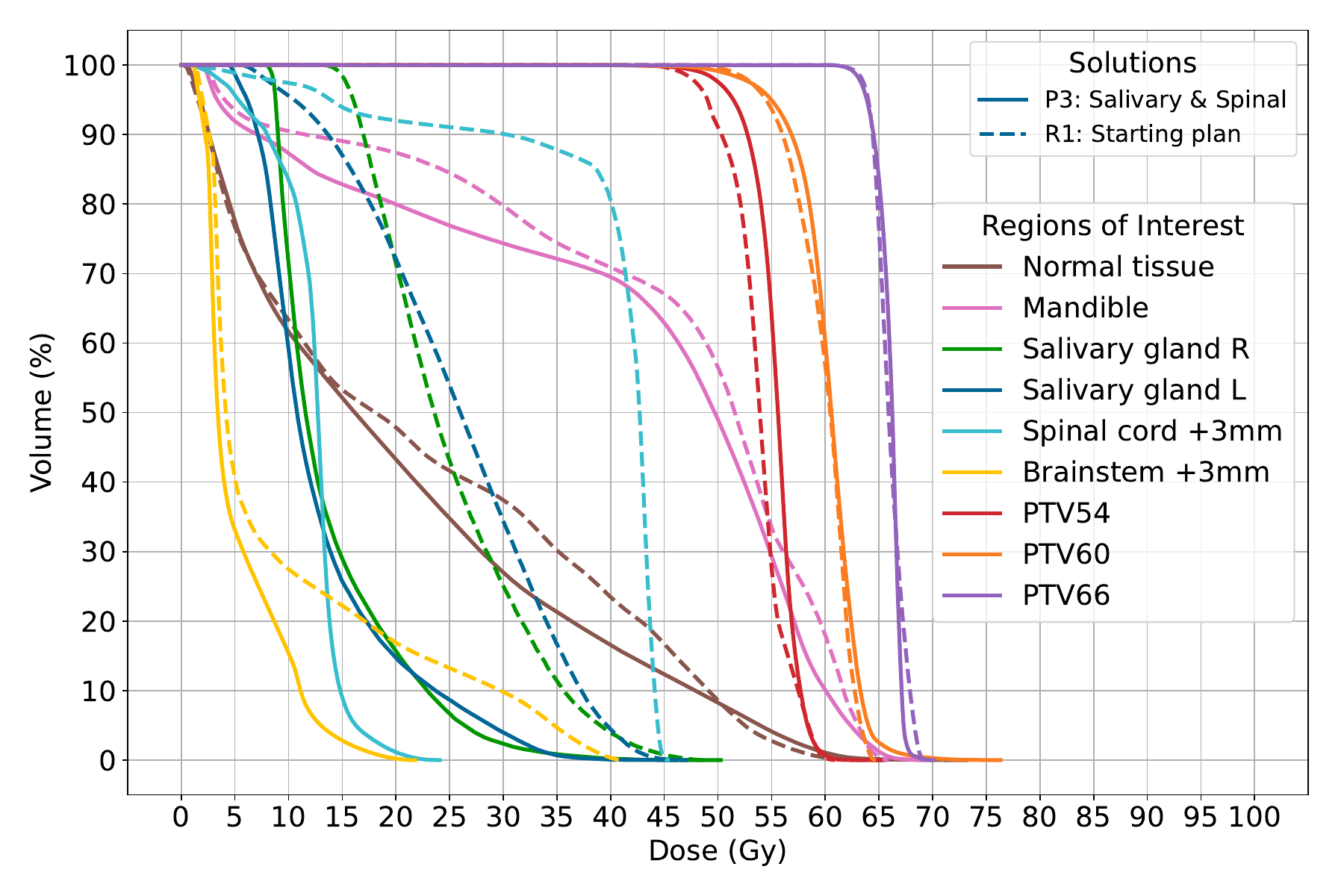}
    \caption{DVHs of Plan 3 (solid) compared with R1, the starting plan (dashed).}
   \label{fig:c3-r1-dvh}
\end{figure}

\FloatBarrier

All plans have to satisfy the same constraints (specified in Table \ref{t:plan3_results}). The comparisons are by DVH and, for illustration, by one representative cross section of the patient model.

The first comparison (Figure \ref{fig:c3-r1-dvh}) is between P3 and R1. As displayed, P3 reduces the doses delivered to all OARs, with comparable or better dose distributions in PTVs.

\begin{figure}[ht]
    \centering
    \includegraphics[width=\textwidth]{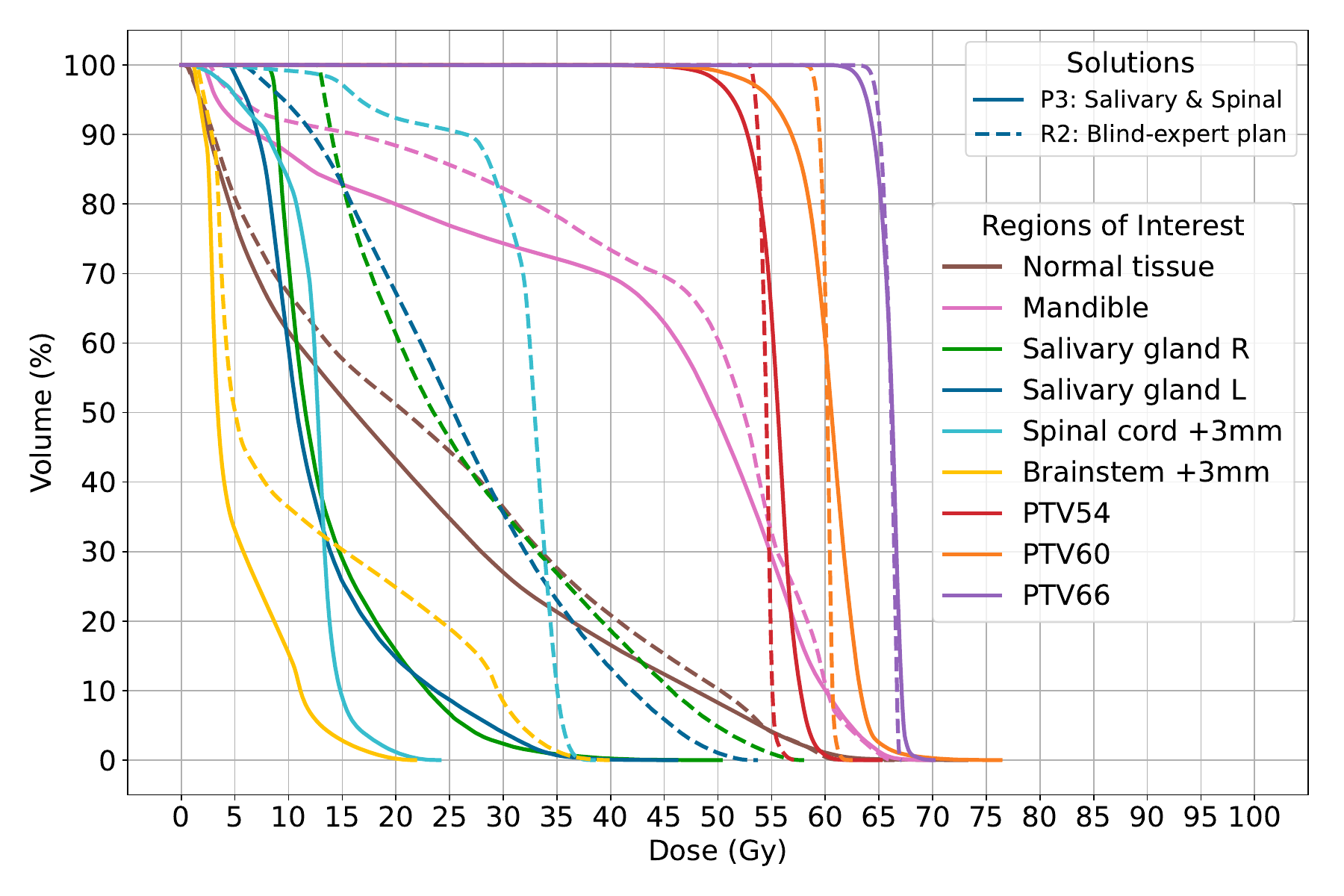}
    \caption{DVHs of plan 3 (solid) compared with the plan R2, the blind-expert plan (dashed).}  
   \label{fig:c3-r2-dvh}
\end{figure}

When comparing P3 with R2 (Figure \ref{fig:c3-r2-dvh}), R2 provides a better dose distribution in PTVs, however at the expense of higher dose deposited in OARs, especially in both salivary glands and the spinal cord +3mm. Notably, in P3 the doses deposited in the salivary glands are 20$\%$ lower than in R2. And in P3 the dose deposited in PTV66 is within the specified bounds. 

\begin{figure}[ht]
    \centering
    \includegraphics[width=\textwidth]{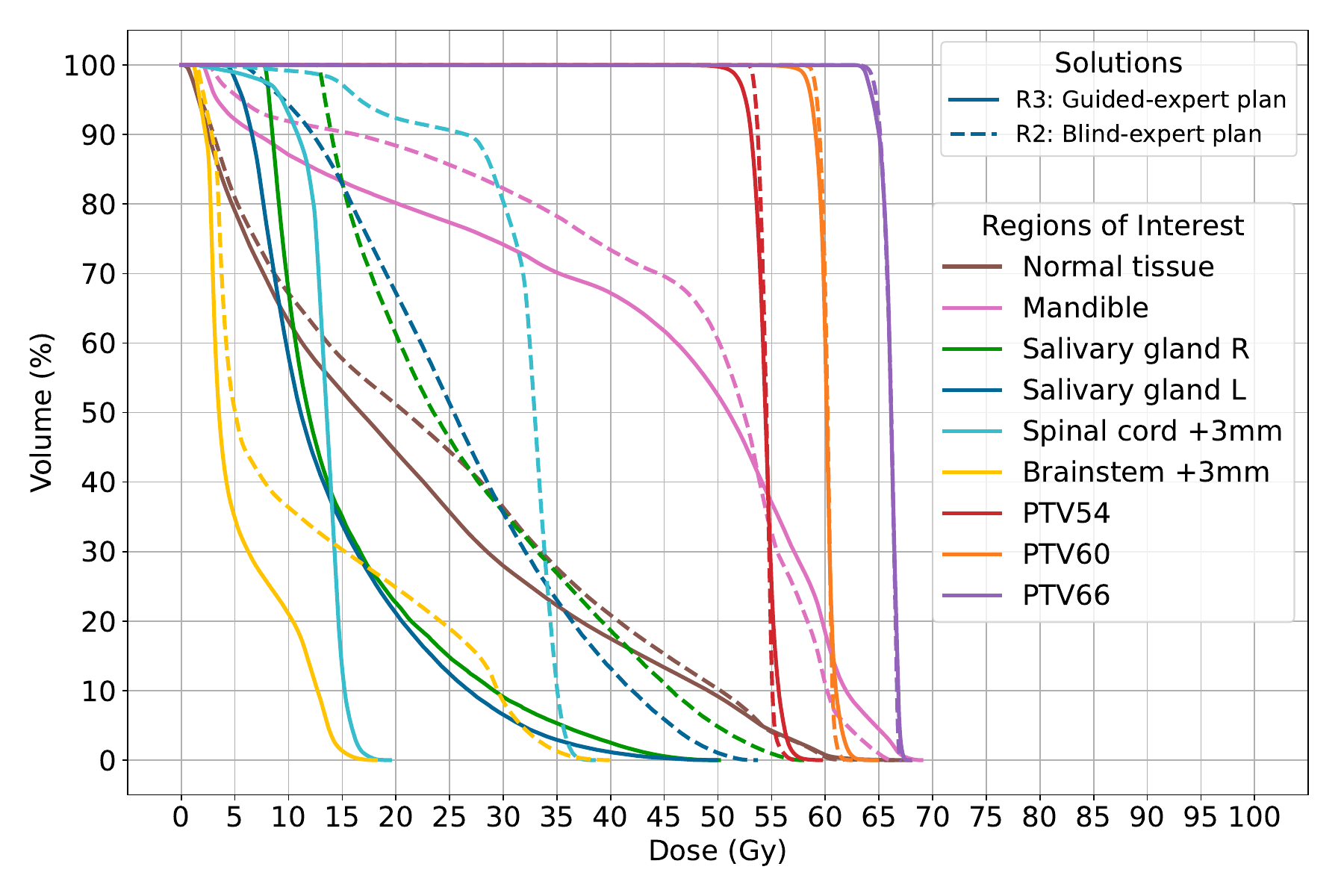}
    \caption{DVHs of the guided plan R3 (solid) compared with R2 (dashed).}  
   \label{fig:r3-r2-dvh}
\end{figure}

\FloatBarrier

Next, we compare with \textit{Guided-expert plan} (R3), a plan prepared by the Expert using P3 as the Eclipse starting plan, with R2 (Figure \ref{fig:r3-r2-dvh}). R3 provides a similar to R2 dose distribution in PTVs and a significant reduction in the doses deposited in OARs. This shows that plans produced by the PersEUD system can be used to enhance RT plans produced by commercial systems currently in use. 

As preliminary conclusions from our work, it can be assumed that using plans from PersEUD as starting points in RT planning systems routinely used in an oncology hospital should improve RT plan quality and reduce RT plan preparation times. This claim has to be validated in more clinical cases.

\begin{figure}[htbp]
    \centering
    \begin{subfigure}[b]{0.44\textwidth}   
        \centering 
        \includegraphics[width=\textwidth]{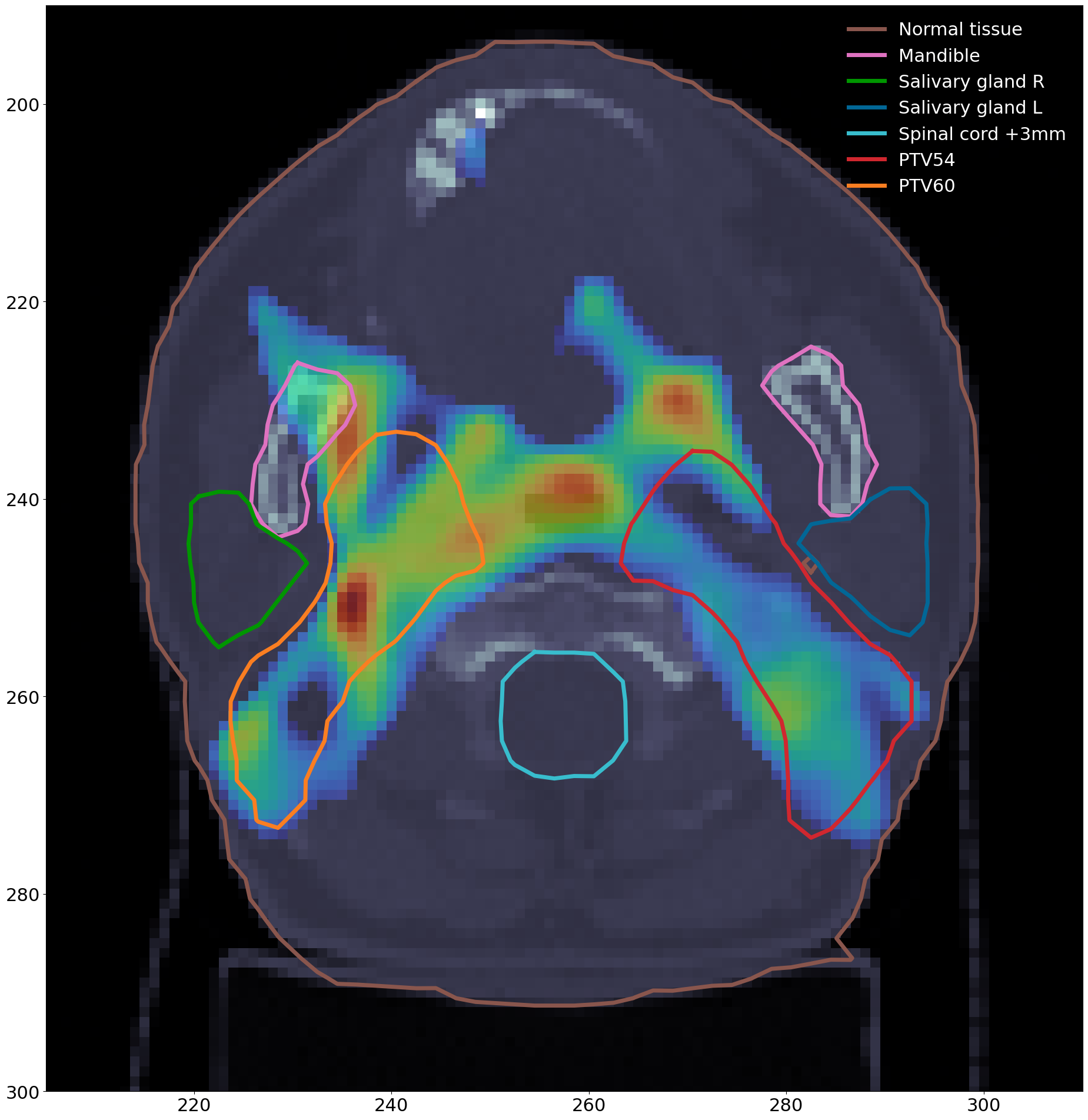}
        \caption{P3: Sal. gl. \& spinal c.}    
        \label{fig:slrk-50}
    \end{subfigure}
    \quad
    \begin{subfigure}[b]{0.44\textwidth}   
        \centering 
        \includegraphics[width=\textwidth]{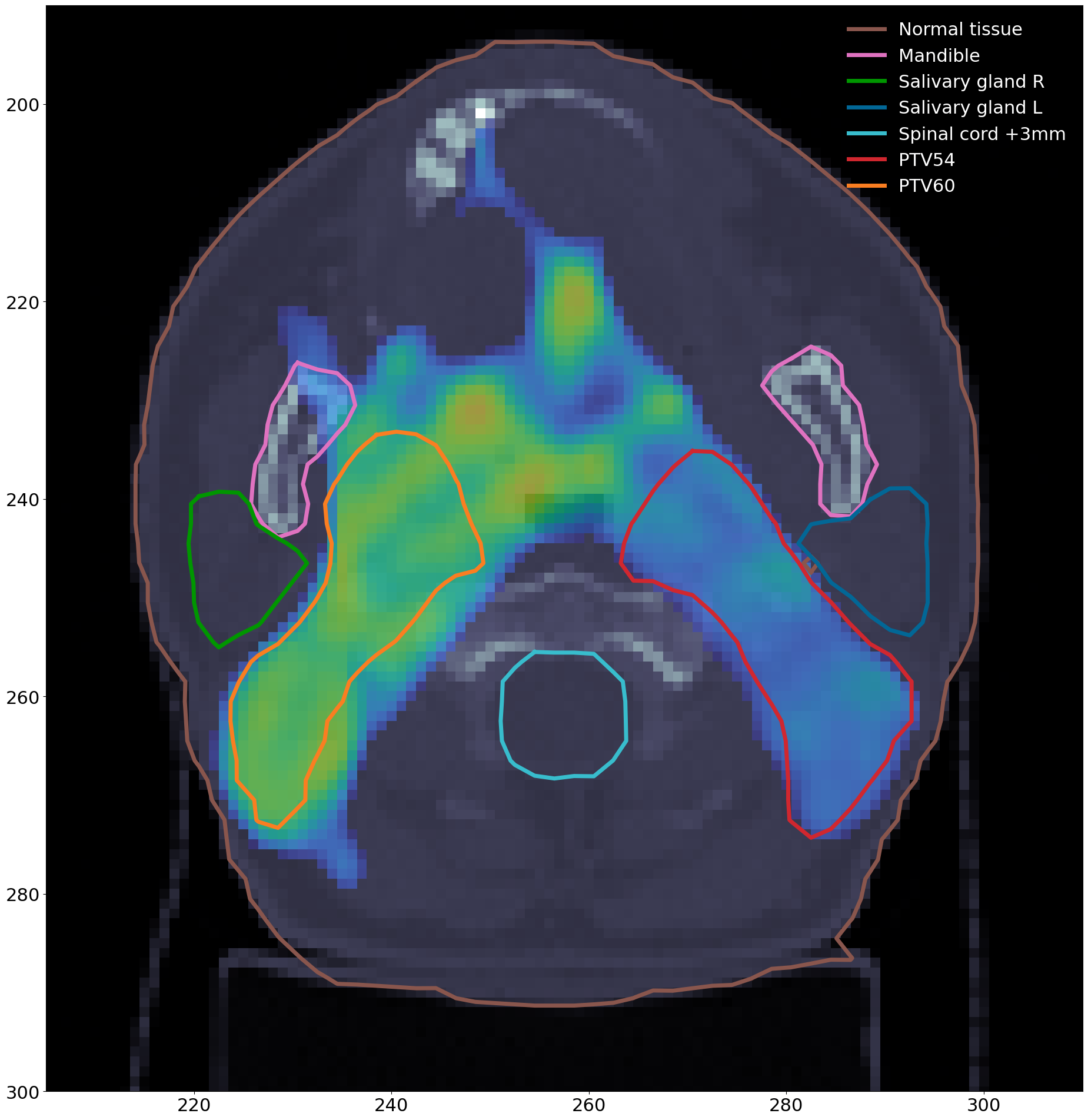}
        \caption{R3: Guided-expert}    
        \label{fig:slrk-re-50}
    \end{subfigure}
    \vskip 0.5\baselineskip
    \begin{subfigure}[b]{0.44\textwidth}   
        \centering 
        \includegraphics[width=\textwidth]{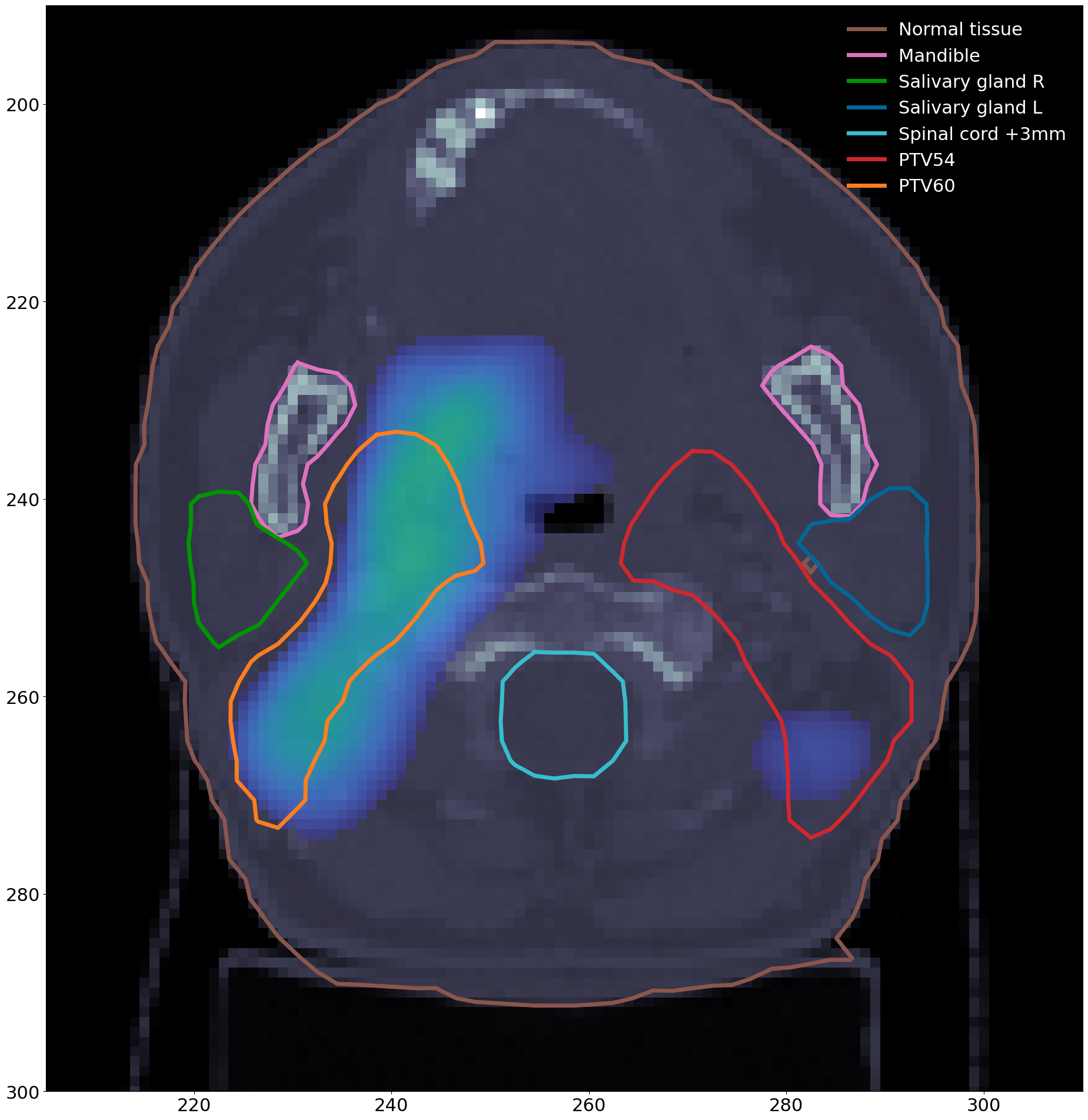}
        \caption{R1: Start}    
        \label{fig:start-50}
    \end{subfigure}
    \quad
    \begin{subfigure}[b]{0.44\textwidth}   
        \centering 
        \includegraphics[width=\textwidth]{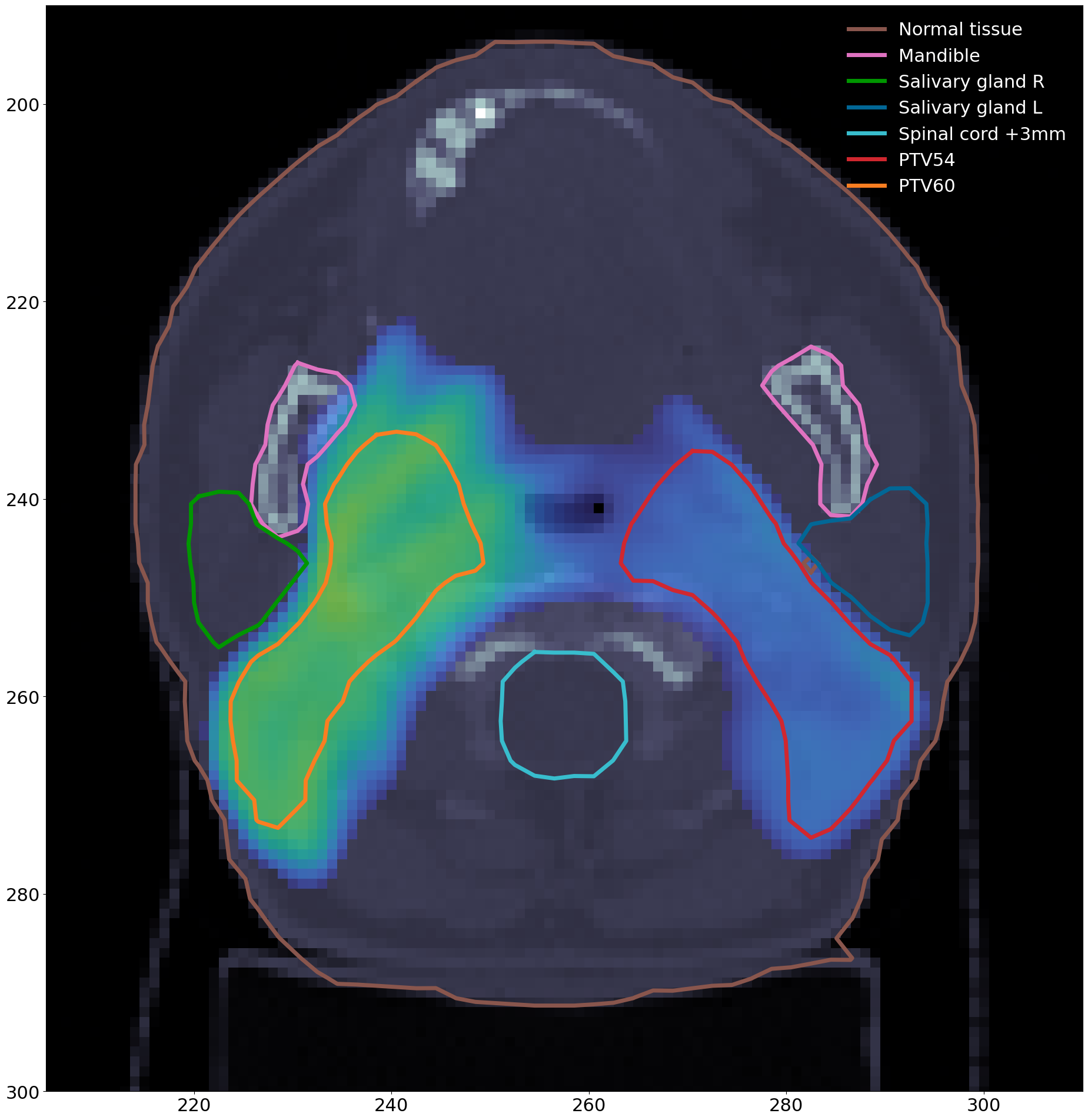}
        \caption{R2: Blind-expert}    
        \label{fig:expert-50}
    \end{subfigure}
    \vskip 0.3\baselineskip
    \begin{subfigure}[b]{\textwidth}   
        \centering 
        \includegraphics[width=0.7\textwidth]{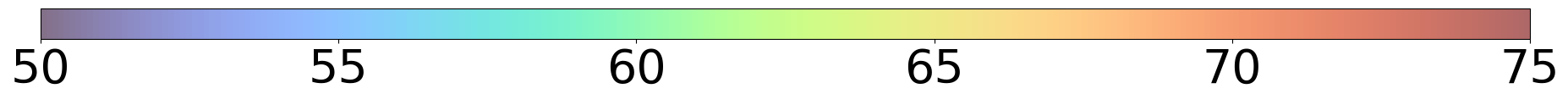}
    \end{subfigure}

    \caption{Dose deposition of the four compared plans for $Z = 83$. Only doses higher than 51.3 Gy are represented.}
    \label{fig:dose-distribution-50}
\end{figure}

\FloatBarrier

\section{Conclusions}

The gEUD-based RT planning optimization has the potential to deliver better RT plans than planning routines based on physical criteria. This is because it directly refers to the biological phenomena in the irradiated cells, whereas the physical criteria measure radiation effects indirectly. Despite its potential, the wide usage of biological optimization techniques in current clinical practices is hampered by the problematic translation between input parameters and resulting dose distributions (RT plans), requiring time-consuming manual model tuning. To help this situation, we propose a system for automated parameter tuning by optimization. To our best knowledge, it is the first system of this sort.

The distinctive future of the proposed system is that it consists of two optimization levels. Because of the number of unknown parameters and the complex nature of the gEUD-based function (\ref{EUDmodel}), this function cannot be optimized directly. So instead, the bi-level optimization scheme has been adopted. On the upper level, the space of the parameters is searched for promising collections of parameters, and the search is guided by values of objective functions of the multi-objective optimization model proposed. Any reasonable heuristic can perform this search, and we have employed evolutionary multi-objective optimization to this aim. On the lower level, the gEUD-based function (function (\ref{EUDmodel})) optimization is performed by a custom-implemented exact optimization method.

The numerical results show the potential benefit of using the proposed system. It seems that it can shorten treatment planning time and improve the quality of treatment plans obtained in commercial treatment planning systems. However, this needs to be tested in more clinical cases, also in other tumor locations.

The bi-level optimization scheme we apply here to RT planning optimization is clearly also applicable to any optimization problem in which the search is needed not only for optimal values of model variables but also for favorable (at best: optimal) model parameters.



\begin{funding}
This work has been supported by Grant PID2021-123278OB-I00 funded by MCIN/AEI/ 10.13039/501100011033 and by ``ERDF A way of making Europe''; and by projects PDC2022-133370-I00 and TED2021-132020B-I00  funded by MCIN/AEI/ 10.13039/5011 00011033 and by European Union Next GenerationEU/PRTR. Sav\'ins Puertas Mart\'in is a fellow of the ``Margarita Salas'' grant (RR\_A\_2021\_21), financed by the European Union (NextGenerationEU).
\end{funding}

\bibliographystyle{infor}
\bibliography{Radiotherapy}

\newpage

\begin{biography}\label{bio1}
\author{J.J. Moreno} is a doctoral researcher in the Department of Informatics at the University of Almería (UAL), Spain. He obtained his B.Sc. and M.Sc. in Computer Science from UAL in 2015 and 2018, respectively. Since 2017, he has been actively engaged in research as a member of the Supercomputing--Algorithms research group at UAL. His areas of research expertise encompass High-Performance Computing (HPC), Radiotherapy optimization, and Electron Tomography image processing.
\end{biography}

\begin{biography}\label{bio2}
\author{S. Puertas-Mart\'in} is a post-doctoral researcher at the Department of Informatics of the University of Almer\'ia, Spain. He is also doing a research stay at the Information School of the University of Sheffield in the United Kingdom. He obtained his Ph.D. in Computer Science at the University of Almer\'ia in 2020. He is a member of the Supercomputing-Algorithms Research Group at that institution. His research interests are Drug Discovery, Global Optimization and High-Performance Computing.
\end{biography}

\begin{biography}\label{bio3}
	\author{J.L. Redondo} holds the position of full professor at the Department of Informatics at the University of Almería, Spain. She earned her Ph.D. in Computer Science from the same university and is an esteemed member of the Supercomputing-Algorithms Research Group. Her research focuses on High-Performance Computing, Global Optimization, and their diverse applications.
\end{biography}

\begin{biography}\label{bio4}
	\author{P.M. Ortigosa} is a full professor of Architecture and Computer Technology of the University of Almer\'ia, Spain. She received  M.Sc degrees in Physics and Electronic Engineering from the University of Granada in 1994 and 1996, respectively, and a Ph.D. in Computer Science from the University of M\'alaga in 1999. She is a member of the Supercomputing-Algorithms Research Group of the University of Almer\'ia. Her research focuses on High-Performance Computing, Metaheuristic Global Optimization, Computational Intelligence, Deep Learning, and the application to several real problems. Recently she has been working on the Internet of Things.
\end{biography}

\begin{biography}\label{bio5}
\author{A. Zawadzka}, a distinguished Medical Physicist, currently serves in the Medical Physics Department at the Maria Skłodowska-Curie National Research Institute of Oncology. In addition to her role, she holds the position of head of the Treatment Planning Laboratory within the same institution. Her academic journey led to the successful completion of a Ph.D. in Medical Physics from the Maria Skłodowska-Curie National Research Institute of Oncology in 2018. A. Zawadzka is actively engaged as a lecturer at Warsaw University and contributes as a Regional Consultant in the field of Medical Physics. Her scientific pursuits are concentrated on treatment planning, dosimetry, and the quality control of radiotherapy treatments.
\end{biography}

\begin{biography}\label{bio6}
\author{P. Kukolowicz} serves as a Medical Physicist at the Medical Physics Department of the Maria Skłodowska-Curie National Research Institute of Oncology, an institution with a rich history dating back to 1934. As a prominent figure, he holds the position of the head of this esteemed institution. Notably, he has been the President of the Polish Society of Medical Physics from 2011 to 2014 and then again from 2014 to 2018. Since 2018, he is a Consultant in Medical Physics for the Ministry of Health. P. Kukolowicz is also recognized for his role as an advisor to 14 Ph.D. candidates in the field of medical physics. His scientific pursuits primarily revolve around treatment planning, dosimetry, and the quality control of radiotherapy treatments.
\end{biography}

\begin{biography}\label{bio7}
\author{R. Szmurło} holds the position of Assistant Professor at the Warsaw University of Technology, Poland. He has actively contributed to over 10 research projects, funded by entities such as the Polish Ministry of Science, commercial enterprises, and European Union grants, all related to computer simulation methods. The subjects of the projects were related to modeling applications of electromagnetic fields in medical treatment, modeling electric impulse power supply systems, artificial intelligence in medicine for Radiotherapy planning, methods of modeling information systems, among others.

\end{biography}

\begin{biography}\label{bio8}
	\author{I. Kaliszewski} is a full professor in the Systems Research Institute of the Polish Academy of Sciences and in Warsaw School of Information Technology. His scientific interests are in optimization, multiple criteria decision making, computer-aided decision making, and also in identification, quantification, and management of risk business organizations.
\end{biography}

 \begin{biography}\label{bio9}
\author{J. Miroforidis} is an assistant professor at the Systems Research Institute of the Polish Academy of Sciences, Poland. He obtained his Ph.D. in Computer Science from that institution in 2010.
His scientific interests are in multi-objective optimization, multiple criteria decision making, and computer-aided
decision making.
\end{biography}

\begin{biography}\label{bio10}
\author{E.M. Garzón} is a full professor at the Department of Informatics of the University of Almer\'ia, Spain. She obtained her Ph.D. in Computer Science from the University of Almer\'ia in 2000. She is the head of the Supercomputing-Algorithms Research Group at that institution. Her research activity is centered on High Performance Computing (HPC) addressed to extend applications of scientific computation. This way her works have been related to different fields and disciplines,

\end{biography}

\end{document}